\documentclass[12pt,a4paper]{article}
\usepackage{amsmath}
\usepackage{amssymb}
\usepackage{amsthm}
\usepackage{amsfonts}
\usepackage{latexsym}
\usepackage[english]{babel}

\theoremstyle{plain}
\newtheorem{thm}{Theorem}[section]
\newtheorem{cor}{Corollary}[section]
\newtheorem{lem}{Lemma}[section]
\newtheorem{prop}{Proposition}[section]
\theoremstyle{definition}
\newtheorem{defn}{Definition}[section]
\newtheorem{exmp}{Example}[section]

\theoremstyle{remark}
\newtheorem{claim}{Claim}[section]
\newtheorem{rem}{Remark}[section]

\title{Moment maps and equivariant volumes}
\author{Alberto Della Vedova and Roberto Paoletti
\footnote{\noindent{\bf Address.} Dipartimento di Matematica e
Applicazioni, Universit\`a degli Studi di Milano Bicocca, Via R.
Cozzi 53, 20125 Milano, Italy; {\bf e-mail}:
alberto.dellavedova@unimib.it, roberto.paoletti@unimib.it }}
\date{}
\begin{document}
\maketitle
\section{Introduction}

Let $M$ an n-dimensional complex projective variety, and consider a
holomorphic line bundle $L$ on $M$. The \textit{volume} of $L$ is
defined as
\begin{equation}
\label{eqn:defn-of-volume} \mathrm{vol}(L)=:\limsup _{k\rightarrow
+\infty } \frac{\mathrm{n}!}{k^\mathrm{n}}\dim H^0(M,L^{\otimes
k}).\end{equation} In the special case where $L$ is ample (or, more
generally, nef and big) $\mathrm{vol}(L)$ is simply the top
self-intersection of $L$, but in general the study of the volume of
arbitrary big line bundles has proven to be a subtle and rich
subject (see \cite{f}, \cite{del}, and \cite{lazI} for a complete
discussion, examples, and references).

Suppose now that $G$ is a compact connected g-dimensional Lie group
acting holomorphically on $M$ (one might equivalently start from a
reductive complex Lie group $\tilde G$, but given the prominent role
of moment maps in our discussion we shall mostly work with a fixed
maximal compact subgroup $G\subseteq \tilde G$). Let us fix a
maximal torus $T\subseteq G$ and a Weyl chamber in its dual Lie
algebra, and let $\Lambda _+=\{\mu\}$ be the corresponding set of
highest weights for $G$. For every highest weight $\mu$, let $V_\mu$
be the finite-dimensional irreducible representation of $G$
associated to $\mu$. If the action of $G$ linearizes to a line
bundle $L$ on $M$, for every integer $k\ge 0$ there is an induced
linear action of $G$ on the space of global holomorphic sections
$H^0(M,L^{\otimes k})$. For every $k\ge 0$, therefore, we have a
$G$-equivariant direct sum decomposition
\begin{equation}\label{eqn:direct-sum-decomposition}
H^0(M,L^{\otimes k})=\bigoplus _\mu H^0(M,L^{\otimes
k})_\mu,\end{equation} where for every $\mu$ the summand $
H^0(M,L^{\otimes k})_\mu$ is $G$-equivariantly isomorphic to a
direct sum of copies of $V_\mu$: $H^0(M,L^{\otimes k})_\mu\cong
V_\mu ^{\oplus N^{(k)}(\mu)}$ as $G$-modules. In the ample case, the
relation between the algebraic structure of the decomposition
(\ref{eqn:direct-sum-decomposition}), expressed by the
multiplicities $N^{(k)}(\mu)$, and the Hamiltonian geometry of the
linearization, expressed by an associated moment map, is the object
of the circle of ideas revolving around the general principle
\textit{quantization commutes with reduction}; we refer to
\cite{gs-gq}, \cite{jk2}, \cite{kir3}, \cite{mein2}, \cite{ms},
\cite{s}, \cite{ggk}. In the algebraic context, the asymptotic
properties of equivariant modules have been studied in \cite{bd},
\cite{b}. For finite group actions, similar questions have been
considered in \cite{pao-ages}.

Now, given that the expected dimension of the quotient of $M$ by $G$
is $\mathrm{ n-g}$, we are led to introduce and study
\textit{equivariant volumes}
\begin{equation}
\label{eqn:defn-of-equivariant-volume} \mathrm{vol}_\mu(L)=:\limsup
_{k\rightarrow +\infty } \frac{(\mathrm{n-g})!}{k^\mathrm{n-g}}\dim
H^0(M,L^{\otimes k})_\mu.\end{equation}
A number of questions are naturally posed: For example, is $
\mathrm{vol}_\mu(L)$ always finite? Is it homogeneous? Can it be
geometrically determined in terms of appropriate top
self-intersections, as in the ample action-free case? Does it only
depend on the numerical equivalence class of the $G$-linearized line
bundle $L$? If so, does it determine a continuous functions on the
$G$-N\'{e}ron-Severi space of $M$, $\mathrm{NS}^G(M) _{\mathbb{R}}$,
and does an estimate similar to the one in Theorem 2.2.44 of
\cite{lazI} still hold? Also, is there an equivariant version of
Fujita approximation for big classes?

The equivariant setting makes these problems not trivial even when
$L$ is ample. In general, for instance, $\mathrm{vol}_\mu(L)$ may
well be infinite for many $\mu$'s - or identically zero for every
$\mu$.

In this article, after determining equivariant volumes for ample and
regular (meaning that the stable and semi-stable loci coincide, and
are non-empty) $G$-linearized line bundles on $M$ by use of the
Riemann-Roch formulae for multiplicities (\cite{gs-gq},
\cite{mein2}, \cite{jk2}), we shall use ideas and techniques from
the algebro-geometric theory of action-free volumes (especially from
\S 2 of \cite{lazI}) and well-known facts from GIT and the theory of
moment maps (\cite{gs-gq}, \cite{dh}, \cite{kir1}, \cite{kir2},
\cite{t}) to give some answers to the previous questions.

\section{Preliminaries.}

Throughout this article, $G$ will denote a g-dimensional compact
connected Lie group, and $\widetilde{G}$ will be its
complexification, a complex reductive Lie group in which $G$ sits a
maximal compact subgroup.

As is well-known, an holomorphic action of $G$ on a complex
projective variety $M$ extends to a holomorphic action of
$\widetilde{ G}$ on $M$, and if the $G$-action linearizes to a line
bundle $L$, then so does the action of $\widetilde{G}$ \cite{gs-gq}.
We shall briefly say that $M$ is a (complex) projective $G$-variety
(or $\widetilde{G}$-variety), and $L$ a $G$-linearized line bundle
on $M$. After \cite{dh}, we shall denote by $ \mathrm{Pic}^G(M)$ the
group of all isomorphisms classes of $G$-linearized line bundles on
$M$.

If $L\in \mathrm{Pic}^G(M)$ is ample and has non-empty semi-stable
locus, we shall generally denote by $M_0(L)$ the GIT quotient
$M/\!/_L \widetilde{G}=M^{\mathrm{ss}}(L)/\widetilde{G}$.

\subsection{Homological and numerical equivalence}

Let $c_1^{(\mathbb{Z})}:\mathrm{Pic}(M)\rightarrow
H^2(M,\mathbb{Z})$ be the map $c_1^{(\mathbb{Z})}\big (
\mathcal{O}_M(D)\big)=:[D]_{\mathrm{hom}}$, and let
$\mathrm{Pic}(M)_0=:\ker (c_1^{(\mathbb{Z})})$. Thus,
$\mathrm{Pic}(M)_0$ parametrizes the line bundles on $M$ that are
deformations of the trivial line bundle.

Next let $c_1^{(\mathbb{R})}:\mathrm{Pic}(M)\rightarrow
H^2(M,\mathbb{R})$ be the composition of $c_1^{(\mathbb{Z})}$ with
the natural map $$ H^2(M,\mathbb{Z})\rightarrow
H^2(M,\mathbb{Z})_{\mathrm{t.f.}}=:H^2(M,\mathbb{Z})/H^2(M,\mathbb{Z})_{\mathrm{tor}}
\hookrightarrow H^2(M,\mathbb{R});$$ here $
H^2(M,\mathbb{Z})_{\mathrm{tor}}$ denotes the torsion part of
$H^2(M,\mathbb{Z})$. Let $\mathrm{Pic}(M)_0'=:\ker
(c_1^{(\mathbb{R})})$. Thus, $L\in \mathrm{Pic}(M)_0'$ if and only
if $c_1^{(\mathbb{Z})}\big ( \mathcal{O}_M(D)\big)\in
H^2(M,\mathbb{Z})_{\mathrm{tor}}$, the torsion part of
$H^2(M,\mathbb{Z})$. Equivalently, $L\in \mathrm{Pic}(M)_0'$ if and
only if some tensor power $L^{\otimes m}$ is a deformation of the
trivial line bundle. Since $H^2(M,\mathbb{Z})_{\mathrm{tor}}$ is
finite, $\mathrm{Pic}(M)_0'$ consists of finitely many connected
components, one for each element of
$H^2(M,\mathbb{Z})_{\mathrm{tor}}$, and every component is
isomorphic to $\mathrm{Pic}(M)_0$.

Given the choice of some Hermitian structure on $L\in
\mathrm{Pic}(M)$, and in view of the de Rahm theorem, we also have
$c_1^{(\mathbb{R})}( L)=\frac{i}{2\pi}[\Theta]$, where $[\Theta]$ is
the cohomology class of the curvature of the compatible connection
on $L$.  Thus, $L\in \mathrm{Pic}(M)_0'$ if and only if $\Theta$ is
an exact form. As shown in 2.3.1 of \cite{dh}, the latter condition
implies that there exists an open trivializing cover for $L$ for
which the transition functions are constants of absolute value $1$.

By the proof of 2.3.3 of \cite{dh}, there is a natural choice of a
$G$-linearization $L^G$ on any $L\in \mathrm{Pic}(M)$ that can be
described by a system of constant transition functions.

Now let $f_G:\mathrm{Pic}^G(M)\rightarrow \mathrm{Pic}(M)$ be the
forgetful map, which to any $G$-linearized line bundle associates
the underlying line bundle (apart from this section, we shall be
systematically sloppy, and use the same symbol for
$L\in\mathrm{Pic}^G(M)$ and its image in $ \mathrm{Pic}(M)$). By the
above considerations, we have a commutative diagram of split short
exact sequences:
\begin{equation}\label{eqn:short-exact-pic}
\begin{array}{ccccccccc}
  1 & \rightarrow & \chi(G) & \longrightarrow & \ker (c_1^{(\mathbb{R})}\circ f_G) & \stackrel{f_G}{\longrightarrow} & \mathrm{Pic}(M)_0' & \rightarrow&0 \\
    &             &         &             &  & &&\\
    &             &   \|    &             &          \uparrow                  &             &  \uparrow          &            & \\
    &&&&&&&\\
  1 & \rightarrow & \chi(G) & \longrightarrow  & \ker (c_1^{(\mathbb{Z})}\circ f_G) & \stackrel{f_G}{\longrightarrow} & \mathrm{Pic}(M)_0  &
  \rightarrow&0,
\end{array}
\end{equation}
where $\chi (G)$ is the character group of $G$, the vertical arrows
are inclusions, and the splitting is given by the maps $L\mapsto
L^G$.

\begin{defn}\label{defn:homolo-trivial-G}
(\cite{dh}, Definition 2.3.4) Let us set
\begin{eqnarray*}
\mathrm{Pic}^G(M)_0&=:&\{ L^G:L\in \mathrm{Pic}(M)_0\},\\
 \mathrm{Pic}^G(M)_0'&=:&\{ L^G:L\in \mathrm{Pic}(M)_0'\}.
 \end{eqnarray*}
Thus $\mathrm{Pic}^G(M)_0$
 is the
connected subgroup of the \textit{homologically trivial}
$G$-linearized line bundles in the sense of \cite{dh}. We shall call
the $G$-linearized line bundles in $\mathrm{Pic}^G(M)_0'$
\textit{numerically trivial}.
\end{defn}

\begin{rem}\label{rem:jacobian-connected}
Clearly $ \mathrm{Pic}^G(M)_0\cong \mathrm{Pic}(M)_0$ and $
\mathrm{Pic}^G(M)_0'\cong \mathrm{Pic}(M)_0'$.\end{rem}

\begin{defn}\label{defn:numerical-vs-homological-equivalence}
We shall say that $L_1,L_2\in \mathrm{Pic}^G(M)$ are
\textit{homologically equivalent}, written $L_1\sim _\mathrm{h}
L_2$, if $L_1\otimes L_2^{-1}\in \mathrm{Pic}^G(M)_0$. We shall say
that $L_1,L_2\in \mathrm{Pic}^G(M)$ are \textit{numerically
equivalent}, written $L_1\sim _\mathrm{n} L_2$, if $L_1\otimes
L_2^{-1}\in \mathrm{Pic}^G(M)_0'$.\end{defn}

Let now $ \mathrm{NS}^G(M)=:\mathrm{Pic}^G(M)/\mathrm{Pic}^G(M)_0$
be the set of all homological equivalence classes (2.3.7 of
\cite{dh}). Then $ \mathrm{NS}^G(M)$ is a finitely generated Abelian
group, of rank $\rho ^G(M)=:\rho (M)+t(\widetilde{G})$, where $\rho
(M)$ is the Picard number of $M$, and $t(\widetilde{G})$ is the
dimension of the radical of $\widetilde{G}$. Let $
\mathrm{NS}^G(M)_{\mathrm{tor}}\subseteq \mathrm{NS}^G(M)$ be the
torsion part of $\mathrm{NS}^G(M)$.

\begin{lem} \label{lem:torsion-is-torsion-pic}
$
\mathrm{NS}^G(M)_{\mathrm{tor}}=\mathrm{Pic}^G(M)_0'/\mathrm{Pic}^G(M)_0$.
\end{lem}

\textit{Proof.} By (\ref{eqn:short-exact-pic}) and Definition
\ref{defn:homolo-trivial-G}, we have isomorphisms of Abelian groups
\begin{eqnarray}\label{eqn:iso-given-splitness}
\ker (c_1^{(\mathbb{R})}\circ f_G)&\cong &\chi (G)\times
\mathrm{Pic}(M)'_0\,\cong \,\chi (G)\times \mathrm{Pic}^G(M)'_0,\\
\ker (c_1^{(\mathbb{Z})}\circ f_G)&\cong &\chi (G)\times
\mathrm{Pic}(M)_0\,\cong \,\chi (G)\times
\mathrm{Pic}^G(M)_0.\nonumber
\end{eqnarray}
Under this isomorphism, for any $L\in \mathrm{Pic}(M)'_0$ its lift
$L^G$ corresponds to the pair $(1,L)$, and for any integer $m$ the
tensor power $(L^G) ^{\otimes m}$ corresponds to the pair
$(1,L^{\otimes m})$. Since $L\in \mathrm{Pic}(M)'_0$, $L^{\otimes m
}\in \mathrm{Pic}(M)_0$, for some integer $m\ge 1$. Therefore,
$(L^G) ^{\otimes m}\in \mathrm{Pic}^G(M)_0$, and the homological
equivalence class $[L]\in \mathrm{NS}^G(M)$ is torsion.

Conversely, suppose $R\in \mathrm{Pic}^G(M)$ and that its class
$[R]$ in $\mathrm{NS}^G(M)$ is torsion. This means that $R^{\otimes
m}\in \mathrm{Pic}^G(M)_0$ for some integer $m\ge 1$. Thus
$f_G(R)^{\otimes m}\in \mathrm{Pic}(M)_0$ and therefore $f_G(R)\in
\mathrm{Pic}(M)'_0$. In other words, $R\in \ker
(c_1^{(\mathbb{R})}\circ f_G)\cong \chi (G)\times
\mathrm{Pic}(M)'_0$.

Suppose then that $R$ corresponds to a pair $(\gamma,L)$, where
$\gamma \in \chi (G)$ and $L\in \mathrm{Pic}_0(M)'$. Then
$L^{\otimes m}$ corresponds to $(\gamma ^m,R^{\otimes m})$. In order
for this to lie in $\mathrm{Pic}^G(M)_0$, we need to have $\gamma
^m=1$, whence $\gamma =1$ (by the connectedness of $G$). Thus, $R$
lies in $ \mathrm{Pic}^G(M)_0'$.

\begin{defn}\label{defn:numerical-lattice} (2.3.7 of \cite{dh})
The $G$-\textit{numerical lattice} of $M$ is
$$\mathrm{Num}^G(M)=:\mathrm{NS}^G(M)/\mathrm{NS}^G(M)_{\mathrm{tor}}.$$
\end{defn}

By Lemma \ref{lem:torsion-is-torsion-pic}, $\mathrm{Num}^G(M)$ is
the set of all numerical equivalence classes:

\begin{cor}\label{cor:num-lattice}
$\mathrm{Num}^G(M)\cong \mathrm{Pic}^G(M)/\mathrm{Pic}^G(M)'_0$.
\end{cor}

\begin{defn} \label{defn:equiv-neron-severi}
The \textit{equivariant N\'{e}ron-Severi space} of $M$ is
$$ \mathrm{NS}^G(M) _{\mathbb{R}}=:\mathrm{NS}^G(M)\otimes
_\mathbb{Z}\mathbb{R}.$$
\end{defn}

\begin{rem}
$ \mathrm{NS}^G(M) _{\mathbb{R}}$ is a real vector space of
dimension $\rho ^G(M)$, it contains the numerical lattice and $
\mathrm{NS}^G(M) _{\mathbb{R}}\cong \mathrm{Num}^G(M)\otimes
_\mathbb{Z}\mathbb{R}$.
\end{rem}

\begin{defn}\label{defn:G-ample-cone}
The \textit{$G$-ample cone} $\mathrm{C}^G(M)\subseteq
\mathrm{NS}^G(M) _{\mathbb{R}}$ is the convex cone spanned by the
classes of all ample $L\in \mathrm{Pic}^G(M)$ having non-empty
semi-stable locus.
\end{defn}

\subsection{Moment maps}\label{subsctn:moment-maps}

In the following, by a \textit{$G$-linearized Hermitian line bundle}
on the projective $G$-variety $M$ we shall mean the choice of $L\in
\mathrm{Pic}^G(M)$ together with a $G$-invariant Hermitian metric
$h$ on $L$. If $M$ is non-singular and $(L,h)$ is a $G$-linearized
Hermitian line bundle on $M$, let $X\subseteq L^*$ be the
corresponding unit circle bundle, and let $\alpha \in \Omega ^1(X)$
be the connection form for the unique compatible connection.  Let
$\mathfrak{g}$ be the Lie algebra of $G$. Then we may define an
equivariant map $\Phi _{L,h}:M\rightarrow \mathfrak{g}^*$ by
setting, for every $m\in M$ and $\xi \in \frak{g}$:
\begin{equation}\label{eqn:moment-map-general}
\left <\Phi _{L,h}(m),\xi\right >=:\alpha _x\big (\xi _X (x)\big ),
\end{equation}
where $x\in X$ is any point lying over $m$, and $\xi _X$ is the
vector field on $X$ generated by $\xi$ \cite{ggk}.

\begin{defn}\label{defn:moment-map-general}
$\Phi _{L,h}$ will be called a moment map for $L$.\end{defn}

When $L$ is ample, $h$ can be so chosen that the curvature form of
the unique compatible connection is $-2\pi i\Omega$, where $\Omega$
is a $G$-invariant K\"{a}hler form. Then the action of $G$ on $M$ is
Hamiltonian for $\Omega$, and $\Phi _{L,h}$ is a moment map in the
ordinary sense. In particular, since any $L\in \mathrm{Pic}^G(M)$
may be factored as $L=A\otimes B^{-1}$, where $A,B\in
\mathrm{Pic}^G(M)$ are ample, we see that for any $\xi \in
\mathfrak{g}$ we have
\begin{equation}\label{eqn:differential-of-moment-map}
d\Phi _{L,h}^\xi =\iota (\xi _M)\,\Omega _{L,h},\end{equation} where
 $\Phi _{L,h}^\xi=:\left <\Phi_{L,h},\xi \right >$, $\xi _M$ denotes
 the vector field on $M$ generated by $\xi$, and $-2\pi i\Omega _{L,h}$ is the curvature
form of the unique compatible connection. In the following, we shall
leave the choice of $h$ implicit, and write $\Phi _L$ for $\Phi
_{L,h}$.

Up to topological obstructions, the information encoded by $\Phi _L$
is equivalent to the assignment of the linearization. More
precisely, the infinitesimal action of $\mathfrak{g}$ on the
sections of $L$ is given by
\begin{equation}\label{eqn:infinitesimal-action}
\xi \cdot \sigma =:\nabla _{\xi _M}\sigma +2\pi i \,\Phi _L^\xi
\cdot \sigma\,\,\,\,\,\,\,\,\,\,(\xi \in \mathfrak{g}, \sigma \in
\mathcal{C}^\infty(M,L)),\end{equation}
where $\nabla$ is the associated covariant derivative (equation
(5.1) of \cite{gs-gq}).

\begin{rem}\label{rem:moment-map-general-sing}
When $M$ is not necessarily non-singular, we may still define a
moment map $\Phi _L=\Phi _{L,h}:M\rightarrow \mathfrak{g}^*$ by
pulling-back $(L,h)$ to some equivariant resolution of singularities
$f:\widetilde{M}\rightarrow M$ \cite{eh}, \cite{ev} and checking
that $\Phi _L$ is well-defined (and smooth) on $M$. Then
(\ref{eqn:differential-of-moment-map}) is satisfied on the smooth
locus of $M$. If $L$ is ample, one may equivalently proceed as
follows: For every $k\gg 0$, the linear series $\left |L^{\otimes
k}\right |$ determines a $G$-equivariant embedding $\varphi
_{L^{\otimes k}}:M\rightarrow \mathbb{P}H^0(M,L^{\otimes k})^*$. Let
us fix a $G$-invariant Hermitian metric on $ H^0(M,L^{\otimes k})$,
so that the linear action of $G$ on $ H^0(M,L^{\otimes k})$ is
unitary and the associated action of $G$ on $
\mathbb{P}H^0(M,L^{\otimes k})^*$ is Hamiltonian with respect to the
Fubini-Study form $ \Omega_{\mathrm{FS}}$. Let $\Phi
_{\mathbb{P}H^0(M,L^{\otimes k})^*}:\mathbb{P}H^0(M,L^{\otimes
k})^*\rightarrow \mathfrak{g}^*$ be the associated moment map, and
let us define $$\Phi _L=:\frac 1k\,\Phi _{\mathbb{P}H^0(M,L^{\otimes
k})^*}\circ \varphi _{L^{\otimes k}}:M\rightarrow \mathfrak{g}^*.$$
Again (\ref{eqn:differential-of-moment-map}) is satisfied on the
smooth locus of $M$ when $\Omega=\frac 1k\,\varphi _{L^{\otimes
k}}^*\big (\Omega_{\mathrm{FS}}\big) $.\end{rem}

Given (\ref{eqn:infinitesimal-action}), the arguments in \S 5 and \S
6 of \cite{gs-gq} leading to Theorems 5.3 and 6.3 of \textit{loc.
cit.} yield the following

\begin{thm}\label{thm:pull-back}
Let $M$ a complex projective $G$-variety, and suppose $L\in
\mathrm{Pic}^G(M)$. Let $\Phi _L:M\rightarrow \frak{g}^*$ be a
moment map for $L$. Let $\mu$ be a maximal weight for $G$. Suppose
that $\mu \not\in \Phi _L(M)\subseteq \frak{g}^*$. Then
$H^0(M,L)_\mu=0$.
\end{thm}

We notice in passing the following:

\begin{lem} \label{lem:pull-back-moment-map}
Let $f:N\rightarrow M$ be an equivariant morphism of projective
$G$-varieties. Suppose that $L\in \mathrm{Pic}^G(M)$ and let $\Phi
_L:M\rightarrow \frak{g}^*$ be a moment map for $L$. Then $\Phi
_L\circ f:N\rightarrow \frak{g}^*$ is a moment map for $f^*(L)\in
\mathrm{Pic}^G(N)$.\end{lem}

\subsection{$G$-invariant effectiveness.}

\begin{defn}\label{defn:G-eff}
Let $M$ a complex projective $G$-variety.

\begin{description}
  \item[i):] A $G$-\textit{invariantly effective divisor}
will be an effective divisor $D$ all of whose irreducible components
are $G$-invariant; we shall write $D\ge_G 0$. In other words,
$D\ge_G0$ means that $D=\sum _ia_iD_i$ where $a_i\ge 0$ and each
$D_i$ is irreducible and satisfies $g\cdot D_i=D_i$ $\forall$ $g\in
G$.

  \item[ii):] A \textit{$G$-invariantly effective} line bundle on $M$
  is an $L\in \mathrm{Pic}^G(M)$
  possessing a non-zero $G$-invariant section $0\neq \sigma \in H^0(M,L)^G$.
\end{description}

\end{defn}

\begin{rem} The notion of $G$-invariant effectiveness for $G$-linearized
line bundles differs from the notion of $G$-effectiveness introduced
in \cite{dh} (which means that $M^\mathrm{ss}(L)\neq
\emptyset$).\end{rem}

\begin{rem}
If $L\in \mathrm{Pic}^G(M)$ is $G$-invariantly effective and $0\neq
\sigma \in H^0(M,L)^G$, then $L=\mathcal{O}_M(Z_\sigma)$ for a
$G$-invariantly effective Cartier divisor $Z_\sigma
=\mathrm{zero}(\sigma)\ge _G0 $, but the converse needn't be true.
Some power $L^{\otimes k}$, $k\ge 1$, is $G$-invariantly effective
if and only if $M^\mathrm{ss}(L)\neq \emptyset$ (that is, if and
only if $L$ is $G$-effective in the sense of \cite{dh}). \end{rem}

\begin{rem}\label{rem:G-eff-morphism}
Suppose that $L\in \mathrm{Pic}^G(M)$ is $G$-invariantly effective,
and pick a non-zero $\sigma \in H^0(M,L)^G$. Given any $H\in
\mathrm{Pic}^G(M)$, tensor product by $\sigma ^{\otimes k}$ defines
a $G$-equivariant injective morphism of line bundles $H^{\otimes
k}\rightarrow \big (H\otimes L\big ) ^{\otimes k}$. Thus for every
maximal weight there is an induced injective linear map $H^0(M,
H^{\otimes k})_\mu\rightarrow H^0(M,\big (H\otimes L\big ) ^{\otimes
k})_\mu$. Passing to $\limsup$, we obtain
$$
\mathrm{vol}_\mu (H)\le \mathrm{vol}_\mu(H\otimes L)$$ for every
maximal weight $\mu$.
\end{rem}

\subsection{Regular $G$-linearized line bundles}

\begin{defn}
\label{defn:regular-l-b} Suppose that $B\in \mathrm{Pic}^G(M)$ is
ample. We shall say that $B$ is \textit{regular} if its stable and
semi-stable loci coincide, and are non-empty:
$M^{\mathrm{s}}(B)=M^{\mathrm{ss}}(B)\neq \emptyset$.
\end{defn}

\begin{rem}\label{rem:regular-open-chamber}
Let $\mathrm{NS}^G(M) _{\mathbb{R}}$ be the equivariant real
N\'{e}ron-Severi space of $M$, and let $\mathrm{C}^G(M)\subseteq
\mathrm{NS}^G(M) _{\mathbb{R}}$ be the $G$-ample cone \cite{dh}.
Then $B\in \mathrm{Pic}^G(M)$ is regular if and only if its
numerical equivalence class $[B]\in \mathrm{NS}^G(M) _{\mathbb{R}}$
lies in the interior of some open chamber of $\mathrm{C}^G(M)$.
\end{rem}

\begin{rem}\label{rem:moment-map}
Simplectically, regularity may be formulated as follows: Let $\Phi
_B:M\rightarrow \frak{g}^*$ be a moment map associated to the
$G$-line bundle $B$. Then $B$ is regular if and only if $0\in
\frak{g}^*$ is a regular value of $\Phi _B$, and $\Phi _B
^{-1}(0)\neq \emptyset$ \cite{kir1}.\end{rem}

\begin{rem}\label{rem:reg-quotient}
If $B\in \mathrm{Pic}^G(M)$ is ample and regular, the GIT quotient
$M_0(B)$ with respect to $B$ is a complex orbifold, and $B$ descends
to a line orbibundle $B_0$ on $M_0(B)$. If $\Omega _B$ is a
$G$-invariant K\"{a}hler form representing $c_1(B)$, it descends to
a K\"{a}hler form $\Omega _0$ on $M_0(B)$ representing $c_1(B_0)$.
\end{rem}

\begin{defn}\label{defn:kirwan-resol}
Let $M$ be a complex projective $G$-manifold, and suppose $L\in
\mathrm{Pic}^G(M)$ is ample and has non-empty stable locus:
$M^\mathrm{s}(L)\neq \emptyset$. Then there exist a $G$-equivariant
birational morphism $f:\widetilde{M}\rightarrow M$ from a complex
projective $G$-manifold $\widetilde{M}$, and an effective
$G$-invariant exceptional divisor $E\subseteq \widetilde{M}$, such
that $f^*(L)^{\otimes k}(-E)\in \mathrm{Pic}^G(\widetilde{M})$ is
ample and regular for every $k\gg 0$ \cite{kir2} . Furthermore, $f$
is an isomorphism over the stable points of $L$. We shall refer to a
morphism with these properties as a \textit{Kirwan resolution} of
$(M,L)$.
\end{defn}

\begin{rem}\label{rem:kirwan-resol}
The $G$-action on $\mathcal{O}_{\widetilde{M}}(E)$ constructed in \S
3 of \cite{kir2} is such that $\mathcal{O}_{\tilde M}(E)$ is
$G$-invariantly effective, i.e. the natural morphism (well-defined
up to a non-zero scalar factor)
$\mathcal{O}_{\widetilde{M}}\rightarrow
\mathcal{O}_{\widetilde{M}}(E)$ is $G$-equivariant (the action of
$G$ on the restriction of $\mathcal{O}_{\widetilde{M}}(E)$ to
$\widetilde{M}\setminus E$ is the product of the action of $G$ on
$\widetilde{M}\setminus E$ with the trivial action on $\mathbb{C}$).
\end{rem}

\subsection{Associated characters.}
\label{subsctn:ass0c-character}

Let us first recall some results and terminology from Appendix B of
\cite{ggk}:

\begin{defn}\label{defn:principal-type}
Suppose that $M$ is a complex projective $G$-manifold. Let
$\mathrm{Conj}(G)_{\mathrm{subgr}}$ be the collection of all
conjugacy classes of subgroups of $G$ (that is, an element of
$\mathrm{Conj}(G)_{\mathrm{subgr}}$ is the collection of all
subgroups of $G$ conjugate to a given one). There exists a subgroup
$K\subseteq G$ such that the stabilizer of a general $p\in M$ is
conjugate to $K$. Then its conjugacy class $(K)\in
\mathrm{Conj}(G)_{\mathrm{subgr}}$ is called the \textit{principal
orbit type} of the action of $G$ on $M$.
\end{defn}

\begin{rem}\label{rem:generic-stabilizer-non-empty-stable}
Suppose that $M$ is a complex projective $G$-manifold, with
principal orbit type $(K)$. If there exists an ample $L\in
\mathrm{Pic}^G(M)$ with non-empty stable locus, then $K\subseteq G$
is a finite subgroup. \end{rem}

Let us momentarily suppose that the stabilizer of a general $p\in M$
is a central subgroup, so that the principal orbit type of the
action consists of a single subgroup $K\subseteq G$.

\begin{defn}\label{defn:associated-character-1}
For every $L\in \mathrm{Pic}^G(M)$, let us denote by $\chi
_{K,L}:K\rightarrow S^1$ the character corresponding to the linear
action of $K$ on the fiber of $L$ at a general $p\in M$.\end{defn}

\begin{rem} \label{rem:character-stabilizer}
Let $K^*=:\mathrm{Hom}(K,S^1)$. Then $\chi:L\in
\mathrm{Pic}^G(M)\mapsto \chi _{K,L}\in K^{*}$ is a group
homomorphism. Obviously, $\ker (\chi)$ contains every
$G$-invariantly effective $L\in \mathrm{Pic}^G(M)$. In particular,
if $L,L'\in \mathrm{Pic}^G(M)$, and there exists a $G$-equivariant
morphism $L\rightarrow L'$, then $\chi_{K,L}=\chi_{K,L'}$. In the
situation of Remark \ref{rem:kirwan-resol}, one has $\chi
_{K,\mathcal{O}_{\tilde M }(E)}=1$.
\end{rem}

\begin{rem} \label{rem:condtion-independ-of-p-chi}
If $K$ is not necessarily central, we can still define $\chi
_{K_p,L}:K_p\rightarrow S^1$ for the general $p\in M$ ($K_p$ is the
stabilizer of $p$). The condition $\chi _{K_p,L}=1$ is independent
of $p$.\end{rem}

\begin{defn}\label{defn:associated-character-2}
Suppose again that the stabilizer of the general $p\in M$ is a
central subgroup $K\subseteq G$. Choose a dominant weight $\mu$. If
$T\subseteq G$ is the chosen maximal torus, let $\mu
_{T}:T\rightarrow S^1$ be the character induced by $\mu$ by
exponentiation, $\mu _{T}:\exp _{T}(\xi)\mapsto e^{2\pi
i<\mu,\xi>}$. Given that $K\subseteq T$, we may define $\mu
_K=:\left .\mu _T\right |_K:K\rightarrow S^1$.
\end{defn}

\begin{rem}\label{rem:altern-descr-of-muK}
An alternative description of $\mu_K$ is as follows: Let
$\mathcal{O}_\mu\subseteq \frak{g}^*$ be the coadjoint orbit of
$\mu$. Since $\mu$ is an integral weight, the natural K\"{a}hler
structure on $\mathcal{O}_\mu$ is in fact a Hodge form, that is, it
represents an integral cohomology class. The associated ample
holomorphic line bundle $A_\mu \rightarrow \mathcal{O}_\mu$ is a
$G$-line bundle in a natural manner. Since the action of $K$ is
trivial on $\mathcal{O}_\mu$, as in Definition
\ref{defn:associated-character-1} the linearization induces the
character $\mu _K$ on $K$.
\end{rem}

\begin{defn}\label{defn:numerically-compatible}
In the hypothesis and notation of Definitions
\ref{defn:associated-character-1} and
\ref{defn:associated-character-2}, we shall say that $L\in
\mathrm{Pic}^G(M)$ and a given dominant weight $\mu$ are
\textit{numerically compatible} if for some $r\in \{1,\ldots,|K|\}$
we have
  $\chi _{K,L}^r \cdot \overline{\mu }_K=1$ (the constant character equal to $1$).
\end{defn}

Obviously, any $L\in \mathrm{Pic}^G(M)$ is numerically compatible
with $\mu =0$ (corresponding to the trivial representation).

\subsection{Equivariant exponents}

\begin{defn} \label{defn:G-exponent}
If $L\in \mathrm{Pic}^G(M)$, the \textit{$G$-semigroup} of $G$ is
$$
\mathbb{N}_G(L)=:\big \{m\ge 0:H^0(M,L^{\otimes m})^G\neq \{0\}\big
\}.$$ Assuming $\mathbb{N}_G(L)\neq \{0\}$, the
\textit{$G$-exponent} of $L$ is the greatest common divisor of the
elements of $\mathbb{N}_G(L)$. Thus $e_G(L)$ divides every element
of $\mathbb{N}_G(L)$, and every sufficiently large integral multiple
of $e_G(L)$ belongs to $\mathbb{N}_G(L)$.

For every maximal weight $\mu$, we shall let
$$ \mathbb{N}_\mu(L)=:\big \{m\in \mathbb{N}:H^0(M,L^{\otimes m})_\mu
\neq \{0\}\big \}.$$
\end{defn}

The semigroup $\mathbb{N}_G(L)$ acts on $ \mathbb{N}_\mu(L)$ by
translations.

\begin{exmp} \label{exmp:associated-char-and-exponent}
The following is a consequence of Theorem
\ref{thm:main-regular-case} below. Suppose that the stabilizer
$K\subseteq G$ of a general $p\in M$ is a central subgroup, and that
$L\in \mathrm{Pic}^G(M)$ is ample and regular. Let $\chi
_{K,L}:K\rightarrow S^1$ be the character associated to the
linearization $L$. Then $e_G(L)$ is the period of $\chi _{K,L}$. If
$\mu$ is a maximal weight, it determines a character $\mu
_K:K\rightarrow S^1$. Then $ \mathbb{N}_\mu(L)\neq \{0\}$ if and
only if there exists $r\in \mathbb{N}$ such that $\chi _{K,L}^r\cdot
\overline \mu_K=1$. In this case, if $k\gg 0$ then $k\in
\mathbb{N}_\mu(L)$ if and only if $k\in r+\mathbb{N}_G(L)$. Thus, up
to a finite number of terms, $\mathbb{N}_G(L)$ consists of all
multiples of the period of $\chi_{K,L}$, and $ \mathbb{N}_\mu(L)$
consists of their translates by $r$.\end{exmp}

\begin{lem}\label{lem:e-G-L-p}
For any $L\in \mathrm{Pic}^G(M)$ and any $p\in \mathbb{N}$, we have
$$
e_G(L^{\otimes p})=\frac{e_G(L)}{\mathrm{gcd}(p,e_G(L))}.$$
\end{lem}

\begin{exmp}\label{exmp:exponent-power}
Referring to the situation of Example
\ref{exmp:associated-char-and-exponent}, $e_G(L^{\otimes p})$ is the
period of $\chi _{K,L^{\otimes p}}$, $|\chi _{K,L^{\otimes p}}|$.
Thus,
\begin{equation*}
e_G(L^{\otimes p})=|\chi _{K,L^{\otimes p}}|=|(\chi
_{K,L})^p|=\frac{|\chi _{K,L}|}{\mathrm{gcd}(p,|\chi
_{K,L}|)}=\frac{e_G(L)}{\mathrm{gcd}(p,e_G(L))}.
\end{equation*}
\end{exmp}

\textit{Proof of Lemma \ref{lem:e-G-L-p}.} Let us set
$e'=:e_G(L)/\mathrm{gcd}\big(p,e_G(L) \big)$; hence $e_G(L)=e'\cdot
\mathrm{gcd}\big(p,e_G(L) \big)$, $p=p'\cdot
\mathrm{gcd}\big(p,e_G(L) \big)$ where $e'$ and $p'$ are relatively
prime. By definition, $H^0(M,L^{\otimes p\, e_G(L^{\otimes
p})\,m})^G\neq \{0\}$ for every $m\gg 0$. Thus,
$$e_G(L)\,|\,p\, e_G(L^{\otimes p})\,\Rightarrow
\,e'\,|\,p'\,e_G(L^{\otimes p})\,\Rightarrow \,e'\,|\,e_G(L^{\otimes
p}).$$

On the other hand, since $p\,e'=e_G(L)\big
(p/\mathrm{gcd}\big(p,e_G(L) \big)$ is an integral multiple of
$e_G(L)$, we have $H^0(M,L^{\otimes pe'm})^G\neq \{0\}$ for every
$m\gg 0$. Therefore, $e_G(L^{\otimes p})\,|\,e'$.

\section{Equivariant volumes and GIT quotients}\label{sctn:regular-case}

In the action-free case, the volume of a nef and big line bundle is
computed by the top self-intersection of the first Chern class. With
this in mind, we shall now examine some special cases where the
equivariant volumes admit a similar interpretation, in terms of
suitable top self-intersections on appropriate GIT quotients of $M$.

For simplicity and ease of exposition, we shall focus on the special
case where the stabilizer $K\subseteq G$ of a general $p\in
\Phi^{-1}(0)$ is a central subgroup $K\subseteq G$, which in the
hypothesis of the present section is necessarily finite. In view of
coisotropic embedding Theorem \cite{gs}, $K$ is then the stabilizer
of a general $p\in M$. The proofs below can however be extended with
no conceptual difficulty to the case of an arbitrary principal orbit
type. This will involve singling out for each maximal weight $\mu$
the kernel $K_\mu\subseteq K$ of the action of $K$ on the coadjoint
orbit of $\mu$, and considering the contribution coming from each
conjugacy class of $K_\mu$.

The first case that we consider is the one of ample and regular
$G$-linearized line bundles (Definition \ref{defn:regular-l-b}).

Let us recall that if $L\in \mathrm{Pic}^G(M)$ is ample and regular,
then the GIT quotient of $M$ with respect to $L$, $M_0(L)$, may be
described as a symplectic reduction with respect to a $G$-invariant
K\"{a}hler form $\Omega$ representing $c_1(L)$. More precisely,
recall that by Remark \ref{rem:moment-map} $0\in \frak{g}^*$ is a
regular value of $\Phi$, and $\Phi ^{-1}(0)\neq \emptyset$.
Furthermore, $P=:\Phi ^{-1}(0)$ is a connected $G$-invariant
codimension $g$ submanifold of $M$, on which $G$ acts locally
freely. There is a natural identification between the GIT quotient
$M_0(L) =M^s(L)/\tilde G$ and the symplectic reduction $P/G$, which
carries an induced structure of K\"{a}hler orbifold. In this manner,
$\Omega$ descends to a K\"{a}hler form (in the orbifold sense)
$\Omega _0$ on $M_0(L)$, representing $c_1(L_0)$ - here $L_0$ is
orbifold line bundle that $L$ descends to on $M_0(L)$.

\begin{thm}
\label{thm:main-regular-case} Let $M$ be a complex projective
manifold, $G$ a compact connected Lie group, $\nu :G\times
M\rightarrow M$ a holomorphic action. Suppose (for simplicity) that
the stabilizer of a general $p\in \Phi^{-1}(0)$ is a central
subgroup $K\subseteq G$. Let $L\in \mathrm{Pic}^G(M)$ be ample and
regular, and let $\Phi =\Phi_L:M\rightarrow \mathfrak{g}^*$ be a
moment map for $L$, associated to the K\"{a}hler form $\Omega$
representing $c_1(L)$. Let $M_0(L)=:\Phi ^{-1}(0)/G$ be the
symplectic quotient with respect to $\Phi_L$. Let $\mu\in \Lambda
_+$ be a dominant weight. If $L$ and $\mu$ are not numerically
compatible (Definition \ref{defn:numerically-compatible}), then
$H^0_\mu(M,L^{\otimes k})=0$ for every $k=1,2,\ldots$. If on the
other hand $L$ and $\mu$ are numerically compatible, then
  $$\mathrm{vol} _\mu(L)=
  \dim (V_\mu)^2\cdot \mathrm{vol}\big (M_0(L),\Omega_0\big )>0.$$
  Here
  $$\mathrm{vol}\big (M_0(L),\Omega_0\big )=:
  \int _{M_0(L)}\Omega_0^{\wedge (n-g)},$$
 \end{thm}

The quickest way to prove this is by applying the Riemann-Roch
formulae for multiplicities due to Meinrenken \cite{mein1},
\cite{mein2}. A fairly elementary analytic approach, based on the
microlocal theory of the Szeg\"{o} kernel from \cite{bs} and
\cite{z}, has been used in \cite{pao-mm} in the case of generically
free actions.

\textit{Proof.} For $k\ge 1$, let us replace $L$ by $L^{\otimes k}$,
the Hodge form $\Omega$ and the moment map $\Phi$ by their multiples
$k\Omega$ and $\Phi _k=:k\Phi$. Given any $\mu \in \frak{g}^*$,
there exists $k_0$ such that $\mu$ is a regular value of $\Phi _k$
if $k\ge k_0$. The relevant asymptotic information about the
multiplicity of $V_\mu$ in $H^0(M,L^{\otimes k})$ may then
determined by computing appropriate Riemann-Roch numbers on the
orbifold obtained as symplectic reduction at the coadjoint orbit
$\mathcal{O}_\mu\subseteq \mathfrak{g}^*$ of $\mu$ \cite{kawa},
\cite{mein2}.

More precisely, since $\mu$ is integral the Kirillov symplectic form
$\sigma _\mu$ on $ \mathcal{O}_\mu$ is a Hodge form. By the Kostant
version of the Borel-Bott theorem, there is an ample line bundle
$A_\mu$ on $\mathcal{O}_\mu$ such that $H^0(\mathcal{O}_\mu, A_\mu)$
is the irreducible representation of $G$ with highest weight $\mu$.

Let $M_\mu ^{(k)}$ be the Weinstein symplectic reduction of $M$ at
$\mu$ with respect to the moment map $\Phi _{k}=k\Phi _L$ ($k\gg
0$). Using the normal form description of the symplectic and
Hamiltonian structure of $(M,\Omega)$ in the neighbourhood of the
coisotropic submanifold $P=\Phi ^{-1}(0)$ \cite{mein2},
\cite{gotay}, \cite{gs}, one can verify that $M_\mu^{(k)}$ is, up to
diffeomorphism, the quotient of $P\times \mathcal{O}_\mu$ by the
product action of $G$. In other words, $M_\mu^{(k)}$ is the fibre
orbibundle on $M_0=:P/G$ associated to the principal $G$-orbibundle
$q:P\rightarrow M_0$ and the $G$-space $\mathcal{O}_\mu$ (endowed
with the opposite K\"{a}hler structure); in particular its
diffeotype is independent of $k$ for $k\gg 0$. Let
$p_\mu:M_\mu^{(k)}\rightarrow M_0$ be the projection.

Let $\theta$ be a connection $1$-form for $q$ (\cite{ggk}, Appendix
B). By the shifting trick, the symplectic structure $\Omega
_\mu^{(k)}$ of the orbifold $M_\mu^{(k)}$ is obtained by descending
the closed 2-form $ k\iota ^*(\Omega)+<\mu,F(\theta)>- \sigma _\mu$
on $P\times \mathcal{O}_\mu$ down to the quotient (the symbols of
projections are omitted for notational simplicity). The minimal
coupling term $<\mu,F(\theta)>- \sigma _\mu$ is the curvature of the
line orbibundle $R_\mu=(P\times \overline A_\mu)/G$ on
$M_\mu^{(k)}$. Thus, $ \Omega _\mu^{(k)}$ is the curvature form of
the line orbibundle $p_\mu^*(L_0^{\otimes k})\otimes R_\mu$.

Let
\begin{eqnarray}\label{eqn:tildeP}
\widetilde{P}_\mu&=:&\{(p,\mu',g)\in P\times \mathcal{O}_\mu\times
G\,: \,
g\cdot (p,\mu')=(p,\mu')\},\nonumber \\
\widetilde{P}_{\mu,K}&=:&P\times \mathcal{O}_\mu\times K.
\end{eqnarray}
Let $G$ act on itself by $h\cdot g=:hgh^{-1}$, and on $
\widetilde{P}_\mu$ by the product action. Since $K$ acts trivially
on $P$ and (being central) on $\mathcal{O}_\mu$, there is a natural
inclusion $\widetilde{P}_{\mu,K}\subseteq \widetilde{P}_\mu$. Now
let $\Sigma_\mu=:\tilde P_\mu/G$, $\Sigma _{\mu,K}=:\tilde
P_{\mu,K}/G=M_\mu^{(k)}\times K$. There is a natural orbifold
complex immersion $\Sigma_\mu\rightarrow M_\mu^{(k)}$, with complex
normal orbi-bundle $N_{\Sigma_\mu}$, and $\Sigma _{\mu,K}\subseteq
\Sigma_\mu$ is the union of the $|K|$ connected components mapping
dominantly (and isomorphically) onto $M_\mu^{(k)}$. The orbifold
multiplicity of $\Sigma _{\mu,K}$ is constant and equal to $|K|$.
Let $L_0$ be the line orbi-bundle on $M_0$ determined by descending
$L$, and let $\tilde L_0$ be its pull-back to $\Sigma _\mu$. Let
$\mathrm{r}$ be the complex dimension of $\mathcal{O}_\mu$, so that
$\dim M_\mu^{(k)}=\mathrm{n-g+r}$. After \cite{mein1} and
\cite{mein2}, the multiplicity $N^{(k)}(\mu)$ of the irreducible
representation $V_\mu$ in $H^0(M,L^{\otimes k})$ is then given by:
\begin{eqnarray}\label{eqn:asymptotic-mult-abelian}
N^{(k)}(\mu)&=&\int _{\Sigma _\mu
}\frac{1}{d_{\Sigma_\mu}}\frac{\mathrm{Td}(\Sigma _\mu)
\mathrm{Ch}^{\Sigma _\mu}(p_\mu^*(L_0^{\otimes k})\otimes R_\mu)}{
\mathrm{D}^{\Sigma _\mu}(N_{\Sigma_\mu} )         }\nonumber\\
&=& \int _{\Sigma _{\mu,K}
}\frac{1}{d_{\Sigma_\mu}}\frac{\mathrm{Td}(\Sigma _\mu)
\mathrm{Ch}^{\Sigma _\mu}(p_\mu^*(L_0^{\otimes k})\otimes R_\mu)}{
\mathrm{D}^{\Sigma _\mu}(N_{\Sigma_\mu} )
}+O(k^{\mathrm{n-g}-1})\nonumber\\
&=&k^{\mathrm{n-g}}\,\frac{1}{|K|}\,\sum _{h\in K}\chi
_{K,L}(h)^k\,\overline{
\mu_K(h)} \int _{M_{\mu}^{(k)}}\,\frac{\big (k\,c_1(L_0) +c_1(R_\mu)\big )}{(\mathrm{n-g}+\mathrm{r})!}^{\mathrm{n-g+r}}\nonumber \\
&&+O(k^{\mathrm{n-g}-1}).
\end{eqnarray}
Now suppose that $\chi _{K,L}^k\cdot \overline{\mu_K}\not \equiv 1$.
Then the action of $K$ on $\jmath^*(L^{\otimes k})\boxtimes
\overline A_\mu$ is not trivial, where $\jmath :P\hookrightarrow M$
is the inclusion. Therefore, the fiber of $p_\mu^*(L_0^{\otimes
k})\otimes R_\mu$ on the smooth locus of $M_0$ is a nontrivial
quotient of $\mathbb{C}$, and $N^{(k)}(\mu)=0$ in this case. If
there exists $k$ such that $\chi _{K,L}^k\cdot
\overline{\mu_K}\equiv 1$, on the other hand, the same condition
holds with $k$ replaced by $k+\ell e$, where $e$ is the period of
$\chi _{K,L}$ and $\ell \in \mathbb{Z}$ is arbitrary. Thus $k$ may
be assumed arbitrarily large. Passing to the original K\"{a}hler
structure of $\mathcal{O}_\mu$ in the computation, and recalling
that $\dim (V_\mu)=(r!)^{-1}\int _{\mathcal{O}_\mu}\,\sigma _\mu
^r$, we easily obtain:
\begin{eqnarray}\label{eqn:asymptotic-mult-abelian}
N^{(k)}(\mu)&=&\dim (V_\mu)\,
\frac{k^{\mathrm{n-g}}}{(\mathrm{n-g})!}\,\int
_{M_0}\,c_1(L_0)^{\wedge (\mathrm{n-g})}+O(k^{\mathrm{n-g}-1}).
\end{eqnarray}

\bigskip

As a Corollary of the proof, we obtain:

\begin{cor}\label{cor:invariant-part-asympt}
In the situation of Theorem \ref{thm:main-regular-case}, let $e=e_L$
be the period of $\chi _{K,L}$. Then for every $l\gg 0$ we have
$$
\dim H^0(M,L^{\otimes el})^G =O(l^{\mathrm{n-g}}).$$
\end{cor}

Let us illustrate Theorem \ref{thm:main-regular-case} by a few
examples.

\begin{exmp}\label{exmp:not-homog}
Consider the linear action $S^1\times \mathbb{C}^2\rightarrow
\mathbb{C}^2$ given by $t\cdot (z_0,z_1)=(tz_0,t^{-1}z_1)$. This
descends to an action on $ \mathbb{P}^1$ with an obvious
linearization to the hyperplane bundle
$\mathcal{O}_{\mathbb{P}^1}(1)$. One easily checks that
$\mathrm{vol}_\mu (\mathcal{O}_{\mathbb{P}^1}(1))=1$ for every $\mu
\in \mathbb{Z}$, but $\mathrm{vol}_\mu
(\mathcal{O}_{\mathbb{P}^1}(2))=0$ whenever $\mu$ is odd. To see
that this agrees with Theorem \ref{thm:main-regular-case}, let us
remark that $K=\{\pm 1\}$, that $\chi
_{K,\mathcal{O}_{\mathbb{P}^1}(1)}:K\rightarrow \mathbb{C}^*$ is the
inclusion, and that for every $\mu \in \mathbb{Z}$ the character
$\mu _K:K\rightarrow \mathbb{C}^*$ is exponentiation by $\mu$. It
follows that $\chi _{K,\mathcal{O}_{\mathbb{P}^1}(1)}\cdot \overline
\mu _K=1$ if $\mu$ is odd, and $\chi
_{K,\mathcal{O}_{\mathbb{P}^1}(1)}^2\cdot \overline \mu _K=1$ if
$\mu$ is even. On the other hand, $\chi
_{K,\mathcal{O}_{\mathbb{P}^1}(2)}=\chi
_{K,\mathcal{O}_{\mathbb{P}^1}(1)}^2=1$, and therefore there does
not exist $r\in \mathbb{N}$ such that $\chi
_{K,\mathcal{O}_{\mathbb{P}^1}(2)}^r\cdot \overline \mu _K=1$ when
$\mu$ is odd.
\end{exmp}

\begin{exmp}\label{exmp:positive-dimensional}
For an example with a positive-dimensional quotient, let us now
start from the linear action $S^1\times \mathbb{C}^3\rightarrow
\mathbb{C}^3$ given by $t\cdot
(z_0,z_1,z_2)=:(t^{-1}z_0,tz_1,tz_2)$. Then with the induced
linearization $ \mathcal{O}_{\mathbb{P}^2}(1)$ is a regular ample
$S^1$-linearized line bundle on $ \mathbb{P}^2$, with associated
moment map $$\Phi
([z_0:z_1:z_2])=\frac{|z_1|^2+|z_2|^2-|z_0|^2}{\|z\|^2}.$$ The map
$[z_0:z_1:z_2]\in \Phi^{-1}(0)\mapsto [z_1:z_2]\in \mathbb{P}^1$ is
well-defined and $S^1$-invariant, and shows that
$\mathbb{P}^2/\!/\mathbb{C^*}=\Phi^{-1}(0)/S^1=\mathbb{P}^1$.
Furthermore, $\mathcal{O}_{\mathbb{P}^2}(1)$ does not descend to a
genuine line bundle on $ \mathbb{P}^1$, but so does $
\mathcal{O}_{\mathbb{P}^2}(2)$.  For every $k=0,1,2,\ldots$ and
$\mu=0,1,2,\ldots$ we have
\begin{equation*}
H^0\big( \mathbb{P}^2,\mathcal{O}_{\mathbb{P}^2}(k)\big )_\mu=\left
\{
\begin{array}{ccc}
  \mathrm{span}\{ z_0^{k-a-b}z_1^az_2^b:a+b=\frac{k+\mu}{2}\}& \mathrm{if} & 2\,|\,k-\mu\ge 0 \\
  0 & \mathrm{otherwise}. &
\end{array}\right.
\end{equation*}
It follows first that setting $\mu =0$, $k=2r$ we obtain
\begin{eqnarray*}\dim H^0(
\mathbb{P}^2,\mathcal{O}_{\mathbb{P}^2}(2r))_0&=&\dim H^0(
\mathbb{P}^2,\mathcal{O}_{\mathbb{P}^2}(2r))^G=\dim
\mathbb{C}[z_1,z_2] _{r} \\
&=&1+r=\dim H^0\big(
\mathbb{P}^1,\mathcal{O}_{\mathbb{P}^1}(r)\big).\end{eqnarray*} This
implies that $\mathcal{O}_{\mathbb{P}^2}(2)$ descends to the
hyperplane line bundle on $ \mathbb{P}^1$, and therefore if we
descend the Fubini-Study form on $\mathbb{P}^2$ to a K\"{a}hler form
$\Omega _0$ on $ \mathbb{P}^1$, the latter satisfies
$\mathrm{vol}(\mathbb{P}^1,\Omega_0)=\int _{\mathbb{P}^1}\Omega
_0=\frac 12$.

Next, if in general $k=\mu+2r$, $r\ge 0$, then
\begin{eqnarray*}
\mathrm{vol}_\mu\left ( \mathcal{O}_{\mathbb{P}^2}(1)\right
)&=&\limsup _{r\rightarrow +\infty} \frac{1}{\mu+2r}\dim H^0\big(
\mathbb{P}^2,\mathcal{O}_{\mathbb{P}^2}(k)\big )_\mu\\
& =&\limsup _{r\rightarrow +\infty} \frac{1}{\mu+2r}\dim
\mathbb{C}[z_1,z_2] _{\mu+r}=\frac 12.
\end{eqnarray*}
\end{exmp}

\begin{exmp} For an non-abelian example,
let us consider the injective group homomorphism
$\alpha:\mathrm{SL}(2)\rightarrow \mathrm{SL}(4)$ given by
$$\alpha (A)=\begin{pmatrix}
  A & 0 \\
  0 & (A^t) ^{-1}
\end{pmatrix}.$$
This defines a linear action of $\mathrm{SL}(2)$ on $\mathbb{C}^4$,
whence an action on $ \mathbb{P}^3$ with a built-in linearization to
the hyperplane line bundle. Let us restrict this to the maximal
compact subgroup $\mathrm{SU}(2)\subseteq \mathrm{SL}(2)$, and
identify $\mathfrak{su}(2)^*\cong \mathfrak{su}(2)$ by means of the
Hermitian form $(A,B)\mapsto \mathrm{trace}(A\overline{B}^t)$. The
moment map $\Phi :\mathbb{P}^3\rightarrow \mathfrak{su}(2)\cong
\mathbb{R}\oplus \mathbb{C}$ is given by
$$\Phi([z_0:z_1:z_2:z_3])=\frac{1}{\|z\|^2}\begin{pmatrix}
  |z_0|^2-|z_1|^2+|z_2|^2-|z_3|^2 \\
  z_1\overline{z_0}+z_2\overline{z_3}
\end{pmatrix}.$$
Arguing in affine coordinates, one can see that $0\in
\mathfrak{su}(2)$ is a regular value of $\Phi$, and clearly $\Phi
^{-1}(0)\neq \emptyset$. Thus, $ \mathcal{O}_{\mathbb{P}^3}(1)$ is a
regular $ \mathrm{SU}(2)$-linearized line bundle on $ \mathbb{P}^3$.
Let us set $V=:\mathbb{C}^2$, $W=:\mathbb{C}^4$. There is an
isomorphism of $ \mathrm{SU}(2)$-modules $W\cong V\oplus V^*\cong
V\oplus \overline V\cong V\oplus V$. Therefore, as $
\mathrm{SU}(2)$-modules,
\begin{eqnarray}\label{eqn:decomp-su2}
H^0\big( \mathbb{P}^3,\mathcal{O}_{\mathbb{P}^3}(k)\big )&\cong&
\mathrm{Sym}^k(V\oplus V)\cong \bigoplus
_{a=0}^k\mathrm{Sym}^a(V)\otimes \mathrm{Sym}^{k-a}(V) \nonumber\\
&\cong &\bigoplus _{a=0}^k\mathrm{Sym}^k(V)\oplus
\mathrm{Sym}^{k-2}(V)\oplus \cdots \oplus
\mathrm{Sym}^{|k-2a|}(V)\nonumber \\
&\cong& \bigoplus _{l=0}^{\left [ \frac k2\right
]}\mathrm{Sym}^{k-2l}(V) ^{\oplus (k-2l+1)}.
\end{eqnarray}
It follows from (\ref{eqn:decomp-su2}) that for every $\mu
=0,1,2,\ldots$ we have
\begin{equation}\label{eqn:decomp-su2-piece}
H^0\big( \mathbb{P}^3,\mathcal{O}_{\mathbb{P}^3}(k)\big )_\mu\,=\,
\left\{
\begin{array}{ccc}
  \mathrm{Sym}^{\mu}(V) ^{\oplus (\mu+1)} & \mathrm{if} & \mu \le k, \,2\,|\,k-\mu\\
  0 & \mathrm{otherwise}. &
\end{array}\right.
\end{equation}
We see from (\ref{eqn:decomp-su2-piece}) that $\dim H^0\big(
\mathbb{P}^3,\mathcal{O}_{\mathbb{P}^3}(k)\big )_\mu =(\mu+1)^2$ if
$0\le \mu\le k$, $2\,|\,k-\mu$. It follows that
$\mathrm{vol}_\mu\big (\mathcal{O}_{\mathbb{P}^3}(1)\big )=(\dim
V_\mu)^2$ for every $\mu=0,1,2,\ldots$. To see that this agrees with
Theorem \ref{thm:main-regular-case}, let us remark that $K=\{\pm
\mathrm{id}_2\}$ and that $\chi
_{K,\mathcal{O}_{\mathbb{P}^3}(1)}(\pm \mathrm{id}_2)=\pm 1$. On the
other hand, for every $\mu =0,1,2,\ldots$ the character $\mu _K$ is
given by $\mu _K(\pm \mathrm{id}_2)=(\pm 1)^\mu$.
\end{exmp}

For future reference, and referring now to the general case (that
is, when the stabilizer of a general $p\in M$ is not necessarily a
central subgroup of $G$), we shall next record a Corollary of the
proof of Theorem \ref{thm:main-regular-case}.

\begin{rem}\label{rem:principal-type}
In the situation of Definition \ref{defn:principal-type}, let us
suppose that there exists $L\in \mathrm{Pic}^G(M)$ ample and regular
(or simply with non-empty stable locus). Then any subgroup
$K\subseteq G$ representing the principal orbit type of the action
is finite. If in addition $L$ is $G$-invariantly effective, that is,
$\dim H^0(M,L)^G>0$, then necessarily $\chi _{K,L}=1$. This
condition is independent of the choice of $K\in (K)$.
\end{rem}

\begin{cor}\label{cor:main-regular-case-but-general}
Suppose that $M$ is a complex projective $G$-manifold, and let
$(K)\in \mathrm{Conj}(G)$ be the principal orbit type of the action
of $G$ on $M$. Let $$ \mathcal{C}\subseteq \mathrm{C}^G(M)\subseteq
\mathrm{NS}^G(M) _{\mathbb{R}}$$ be an open chamber in the $G$-ample
cone (so that every $L\in \mathrm{Pic}^G(M)$ whose numerical class
$[L]$ belongs to $ \mathcal{C}$ is ample and regular). Then there
exists a constant $C=C_{\mathcal{C}}>0$ such that the following
holds: Suppose $L\in \mathrm{Pic}^G(M)$, $[L]\in \mathcal{C}$, and
let $(K)$ be the principal orbit type of the $G$-action on $\Phi
^{-1}_L(0)$. If $\chi _{K,L}=1$, then
\begin{equation}\label{eqn:main-regular-case-but-general}
\dim H^0(M,L^{\otimes
k})^G=C\,\mathrm{vol}\big(M_0(\mathcal{C}),\Omega_0(L)\big)\,
k^{\mathrm{n-g}}+\mathrm{L.O.T.}\end{equation}
as $k\rightarrow +\infty$.
\end{cor}

\begin{rem}\label{rem:main-regular-case-but-general}
Here $ M_0(\mathcal{C})$ denotes the GIT quotient of $M$ with
 respect to any of the GIT-equivalent line bundles whose numerical
 class belongs to the open chamber $\mathcal{C}$.
 In particular, the K\"{a}hler form $
\Omega_0(L)$ on $M_0(\mathcal{C})$ depends linearly on $L$ as $L$
varies among the $G$-linearized line bundles
 class belongs to $\mathcal{C}$. Furthermore, as we argue below
the principal orbit type $(K)$ does not depend on the choice of
$[L]\in \mathcal{C}$, so that the condition $\chi _{K,L^{\otimes
k}}=1$ is satisfied by any $[L]\in \mathcal{C}$, for some fixed
integer $k>0$.
\end{rem}

\textit{Proof of Corollary \ref{cor:main-regular-case-but-general}.}
For $[l]\in \mathcal{C}$, let $\Phi_{[l]}:M\rightarrow
\mathfrak{g}^*$ denote the associated Hamiltonian structure. Let us
set $P_{[l]}=(\Phi_{[l]}) ^{-1}(0)$. Then the submanifolds
$P_{[l]}=\Phi ^{-1}_{[l]}(0)$ form a smooth family of $G$-manifolds
as $[l]$ varies in $\mathcal{C}$. In view of the rigidity of compact
group actions \cite{ggk}, this family is locally equivariantly
trivial. Thus, $(K)$ is the principal orbit type of any $\Phi
^{-1}_{[l]}(0)$.

For $[l]\in \mathcal{C}$, let us now define
$\widetilde{P}_{[l],\mu}$ as in (\ref{eqn:tildeP}), and let
$\widetilde{P}_{[l],\mu,(K)}=\bigcup
_i\widetilde{P}_{[l],\mu,(K),i}$ be the union of the connected
components of $ \widetilde{P}_{[l],\mu}$ mapping dominantly onto
$M_0(L)=M_0( \mathcal{C})$. Let us set again $\Sigma
_{[l],\mu}=\widetilde{P}_{[l],\mu}/G$, $\Sigma
_{[l],\mu,K}=\widetilde{P}_{[l],\mu,(K)}/G$, $\Sigma
_{[l],\mu,(K),i}=:\widetilde{P}_{[l],\mu,(K),i}/G$.

For every $i$, let $\mathrm{d}_i$ the cardinality of the generic
stabilizer subgroup of a point in $\widetilde{P}_{\mu,(K),i}$, and
let $ \mathrm{g}_i$ be the degree of the unramified covering $\Sigma
_{[l],\mu,(K),i}\rightarrow M_0( \mathcal{C})$. By the above, the
$\Sigma _{[l],\mu}$ form a smooth family as $[l]$ varies in
$\mathcal{C}$. It follows that the sum $C=\sum
_i\mathrm{g}_i/\mathrm{d}_i$ only depends on $\mathcal{C}$. The
statement follows from this, by arguing as in the proof of Theorem
\ref{thm:main-regular-case}.

\bigskip

Let us now only assume that $L\in \mathrm{Pic}^G(M)$ is ample and
non-empty stable locus. Let $f:\widetilde{M}\rightarrow M$ be a
Kirwan resolution of $(M,L)$. Thus, $f^*(L)^{\otimes k}(-E)\in
\mathrm{Pic}^G\big (\widetilde{M}\big )$ is ample and regular for
every $k\gg 0$, where $E$ is some $f$-exceptional $G$-invariantly
effective divisor. Let $\widetilde{C}$ be the chamber in the
$G$-ample cone of $\widetilde{M}$ containing the numerical class of
$f^*(L)^{\otimes k}(-E)$, $k\gg 0$ \cite{dh}, \cite{t}. Let
$\widetilde{M}_0(\widetilde{C})$ be GIT quotient of $\widetilde{M}$
with respect to $\widetilde{C}$. The line bundle $f^*(L)\in
\mathrm{Pic}^G\big ( \widetilde{M}\big )$ descends to a line
orbi-bundle $L_0$ on $\widetilde{M}_0(\widetilde{C})$; let $
\widetilde{\Omega}_0$ be a 2-form representing its first Chern
class.

\begin{cor}\label{cor:not-empty-stable-locus-case}
Suppose that $L\in \mathrm{Pic}^G(M)$ is ample and has non-empty
stable locus. Let $\mu$ be a dominant weight. Assume again for
simplicity that the stabilizer of a general $p\in M$ is a central
subgroup. If $L$ and $\mu$ are not numerically compatible, then
$H^0_\mu(M,L^{\otimes k})=0$ for every $k=1,2,\ldots$. If on the
other hand $L$ and $\mu$ are numerically compatible, then
  $$\mathrm{vol} _\mu(L)=
  \dim (V_\mu)^2\cdot
  \mathrm{vol}\big (\widetilde{M}_0(\widetilde{C}),\widetilde{\Omega}_0\big )>0.$$
  \end{cor}

\textit{Proof.} Suppose $p\gg 0$ is such that $L^{\otimes
pe_G(L)}(-E)$ is ample and regular; we may then find $r\gg 0$ such
that $H^0\big (\widetilde{M},L^{\otimes pre_G(L)}(-rE)\big )^G\neq
0$. Let $p_i\uparrow +\infty$ be a sequence of integers prime with
$e_G(L)$. It follows from Proposition \ref{prop:key-homog-prime}
below that $\mathrm{vol}_\mu \big (L^{\otimes
p_i}\big)=p_i^{\mathrm{n-g}}\mathrm{vol}_\mu(L)$. Therefore,
\begin{equation}\label{eqn:kirwan-new-limit}
\frac{\mathrm{vol}_\mu\Big (f^*(L)^{\otimes p_i}(-E)\Big
)}{p_i^{\mathrm{n-g}}}\,\le \,\mathrm{vol}_\mu (L )\,\le
\frac{\mathrm{vol}_\mu\Big (f^*(L)^{\otimes p_i+pre_G(L)}(-r E)\Big
)}{p_i^{\mathrm{n-g}}}.\end{equation} The statement follows from
(\ref{eqn:kirwan-new-limit}) by taking the limit as $i\rightarrow
+\infty$.

\bigskip

Let us now briefly consider the more general case where $L\in
\mathrm{Pic}^G(M)$ is nef and big, and has positive 0-volume. Let us
suppose, for simplicity, that there exists $B\in \mathrm{Pic}^G(M)$
ample and such that $M^\mathrm{s}(B)\neq \emptyset$,
$\mathrm{codim}\big (M^\mathrm{u}(B)\big)\geq 2$. Then by Corollary
\ref{cor:nef-and-big-case} below there exist a $G$-invariantly
effective $E\in \mathrm{Pic}^G(M)$ such that $L^{\otimes k}(-E)$ is
ample and $M^\mathrm{s}\Big (L^{\otimes k}(-E)\Big )\neq \emptyset$
for every $k\gg 0$. After replacing $E$ by some power, we can fix
$k_0\gg 0$ divisible by $e_G(L)$ and such that $H^0\big(M,L^{\otimes
k_0}(-E)\big )^G\neq \{0\}$. Letting $p_i\uparrow +\infty$ be a
sequence of integers prime with $e_G(L)$, we obtain
\begin{equation*}\label{eqn:nef-and-big}
\frac{\mathrm{vol}_\mu\Big (L^{\otimes p_i}(-E)\Big
)}{p_i^{\mathrm{n-g}}}\,\le \,\mathrm{vol}_\mu (L )\,\le
\frac{\mathrm{vol}_\mu\Big (L^{\otimes p_i+k_0}(-E)\Big
)}{p_i^{\mathrm{n-g}}}.\end{equation*} If $B$ is regular, it may be
assumed that so is $L^{\otimes k}(-E)$, and that it eventually
belongs to a fixed chamber in the $G$-ample cone of $M$. Then by
passing to the limit as $i\rightarrow +\infty$ we obtain a statement
formally similar to the one in Theorem \ref{thm:main-regular-case}.
Otherwise, we may still apply a Kirwan resolution and reach a
statement as in Corollary \ref{cor:not-empty-stable-locus-case}; we
shall leave it to the reader to complete the argument.

\section{Finiteness}

Unlike the good cases considered in section \S
\ref{sctn:regular-case}, in general equivariant volumes and their
dependence on the weight $\mu$ will have a less controlled behavior.

\begin{exmp}\label{exmp:trivial}
Suppose $\dim (G)>0$. Let $M$ be a complex projective variety, and
$L$ an ample line bundle on $M$. Let $G$ act trivially on $M$ and
$L$. Then $\mathrm{vol}_0(L)=+\infty$ and $\mathrm{vol}_\mu(L)=0$
for all $\mu \neq 0$. Notice that in this case
$M^{\mathrm{u}}(L)=M^{\mathrm{s}}(L)=\emptyset$, and
$M=M^{\mathrm{ss}}(L)$.\end{exmp}

\begin{exmp}\label{exmp:less-trivial}
Consider the linear action $S^1\times
\mathbb{C}^{\mathrm{n}+1}\rightarrow \mathbb{C}^{\mathrm{n}+1}$
given by $$t\cdot (z_0,\ldots,z_\mathrm{n})=(z_0,t z_1,\ldots,t
z_\mathrm{n}).$$ This defines an action on $\mathbb{P}^{\mathrm{n}}$
with an obvious linearization to the hyperplane bundle. If
$\mathrm{n}>1$ then $ \mathrm{vol}_\mu \left
(\mathcal{O}_{\mathbb{P}^{\mathrm{n}}}(1)\right )=0$ for every
$\mu$. In this case,
$M^{\mathrm{ss}}(L)=\mathbb{P}^{\mathrm{n}}\setminus \{z_0=0\}$,
$M^{\mathrm{s}}(L)=\emptyset$.
\end{exmp}

\begin{exmp}\label{exmp:less-less-trivial}
Consider next the linear action $S^1\times
\mathbb{C}^{\mathrm{n}+1}\rightarrow \mathbb{C}^{\mathrm{n}+1}$
given by $$t\cdot (z_0,\ldots,z_\mathrm{n})=(z_0, z_1,\ldots,
z_{\mathrm{n}-1},tz_\mathrm{n}),$$ with the induced action on $
\mathbb{P}^\mathrm{n}$ and the tautological linearization on the
hyperplane bundle. Now $M^{\mathrm{u}}\big
(\mathcal{O}_{\mathbb{P}^{\mathrm{n}}}(1)\big )=
\{[0:\cdots:0:1]\}$, $M^{\mathrm{s}}\big
(\mathcal{O}_{\mathbb{P}^{\mathrm{n}}}(1)\big )=\emptyset$. We have
$ \mathrm{vol}_\mu \left
(\mathcal{O}_{\mathbb{P}^{\mathrm{n}}}(1)\right )=1$ for every $\mu
\ge 0$, $ \mathrm{vol}_\mu \left
(\mathcal{O}_{\mathbb{P}^{\mathrm{n}}}(1)\right )=0$ if $\mu <0$.
\end{exmp}

We shall now give some general hypothesis guaranteeing that every
$G$-line bundle on $M$ satisfies $\mathrm{vol}_\mu(L)<+\infty$ for
every maximal weight $\mu$.

\begin{prop}\label{prop:reg-implies-finite}
Suppose that there exists a regular ample $G$-line bundle $B$ on
$M$. Then $ \mathrm{vol}_\mu(L)<+\infty$ for every $L\in
\mathrm{Pic}^G(M)$ and every maximal weight $\mu$.
\end{prop}

\textit{Proof.} By Corollary \ref{cor:invariant-part-asympt}, we may
find $k=e l\gg 0$ such that $H^0(M,B^{\otimes k})^G\neq \{0\}$. By
Remark \ref{rem:G-eff-morphism}, $ \mathrm{vol}_\mu (L)\le
\mathrm{vol}_\mu (L\otimes B^{\otimes k})$.

Hence it suffices to show that if $k=el\gg 0$ then $
\mathrm{vol}_\mu (L\otimes B^{\otimes k})<+\infty$ for every $\mu$.

If $k\gg 0$ the $G$-line bundle $L\otimes B^{\otimes k}$ is ample
and its equivalence class lies in the interior of the same chamber
as $B$; in particular, it is regular. The statement follows by
Theorem \ref{thm:main-regular-case}.

\begin{cor}\label{cor:stable-is-enough}
Suppose that there exists an ample $G$-line bundle $A$ on $M$ such
that $M^{\mathrm{s}}(A)\neq \emptyset$. Then $\mathrm{vol}_\mu
(L)<+\infty$ for every $L\in \mathrm{Pic}^G(M)$ and every maximal
weight $\mu$.
\end{cor}

\textit{Proof.} Let $f:\tilde M\rightarrow M$ be a Kirwan resolution
of the pair $(M,A)$ (Definition \ref{defn:kirwan-resol}). By
definition, there exists an ample and regular $G$-line bundle on
$\tilde M$. Hence $\mathrm{vol}_\eta(N)<+\infty$ for every $G$-line
bundle $N$ on $\tilde M$ and every maximal weight $\eta$.

Now suppose $L\in \mathrm{Pic}^G(M)$. There is a $G$-equivariant
isomorphism $f_*\big(f^*(L)\big ) ^{\otimes k}\cong L^{\otimes k}$,
whence $G$-equivariant isomorphisms $$H^0(\tilde M,f^*(L)^{\otimes
k})\cong H^0(M,L^{\otimes k})$$ for every $k=1,2,\ldots$. Therefore,
given any maximal weight $\mu$, we have
$\mathrm{vol}_\mu(L)=\mathrm{vol}_\mu \big(f^*(L)\big)<+\infty.$

\bigskip

If on the other hand an ample $G$-line bundle $L$ has no semi-stable
points, then all equivariant volumes of $L$ vanish.

\begin{prop}\label{prop:unst-and-vanish}
Let $M$ be a projective $G$-variety. Let $L\in \mathrm{Pic}^G(M)$ be
ample and such that $M=M^\mathrm{u}(L)$. Then for every maximal
weight $\mu$ there exists $r_\mu\ge 0$ such that $H^0(M,L^{\otimes
r})_\mu =0$ for every integer $r\ge r_\mu$. In particular,
$\mathrm{vol}_\mu(L)=0$ for every maximal weight $\mu$.\end{prop}

\textit{Proof of Proposition \ref{prop:unst-and-vanish}.} Let $\Phi
_L:M\rightarrow \frak{g}^*$ be a moment map for $L$ (Definition
\ref{defn:moment-map-general} and Remark
\ref{rem:moment-map-general-sing}); then $\Phi _{L^{\otimes
r}}=:r\Phi _L$ is a moment map for $L^{\otimes r}$. By the
hypothesis, $\Phi _L$ is bounded in norm away from zero \cite{kir1}.
Therefore, there exists $c>0$ such that
\begin{equation}\label{eqn:bound-grows-with-r}
\min \{ \|\Phi _{L^{\otimes r}}(m)\|:m\in M\}\ge cr.
\end{equation}
Given any $\mu \in \frak{g}^*$ we then have $\mu \not \in r\,\Phi
_L(M)=\Phi _{L^{\otimes r}}(M)$ if $r\gg 0$.

If $M$ is a projective manifold, Theorem 6.3 of \cite{gs} implies
$H^0(M, L^{\otimes r})_\mu =\{0\}$ if $r\gg 0$.

If on the other hand $M$ is singular, let us choose a
$G$-equivariant resolution of singularities, $f:\tilde M\rightarrow
M$. There exists a $G$-invariantly effective exceptional normal
crossing divisor $E\subseteq \tilde M$ such that $f^*(L) ^{\otimes
r}(-E)\in \mathrm{Pic}^G(\tilde M)$ is ample for every $r\gg 0$. Let
us consider the $G$-equivariant short exact sequence
\[0\longrightarrow f^*(L) ^{\otimes r}(-E)
\longrightarrow f^*(L) ^{\otimes r}\longrightarrow
\mathcal{O}_E\otimes f^*(L) ^{\otimes r}\longrightarrow 0.\] At
every maximal weight $\mu$, this induces an exact sequence
\[0\longrightarrow H^0\big(\tilde M,f^*(L) ^{\otimes r}(-E)\big )_\mu
\longrightarrow H^0(\tilde M,f^*(L) ^{\otimes r})_\mu\longrightarrow
H^0(E,\mathcal{O}_E\otimes f^*(L) ^{\otimes r})_\mu.\]
A moment map for $f^*(L) ^{\otimes r}(-E)$ has the form $\Phi _r=
r\,\Phi _L\circ f+\Phi _0:\tilde M\rightarrow \frak{g}^*$, where
$\Phi _0:M\rightarrow \frak{g}^*$ is a constant equivariant map. It
follows from (\ref{eqn:bound-grows-with-r}) that $\mu \not\in \Phi
_r(\tilde M)$ if $r\gg 0$, and therefore $H^0(\tilde M, f^*(L)
^{\otimes r}(-E))_\mu =0$, again by Theorem 6.3 of \cite{gs}.

Next let us decompose $E$ as the sum of its irreducible components:
$E=\sum _jE_j$. Each $E_j$ is an irreducible $G$-invariant subscheme
of $\tilde M$, and there is an equivariant injection
$\mathcal{O}_E\rightarrow \bigoplus _j\mathcal{O}_{E_j}$. Therefore,
it suffices to prove that $H^0(E_j,\mathcal{O}_{E_j}\otimes f^*(L)
^{\otimes r})_\mu=0$ for every $j$ and every $r\gg 0$. If all the
$E_j$'s are reduced, hence non-singular varieties, the statement
follows from Theorem \ref{thm:pull-back}. If some of the $E_j$'s are
not reduced, the same conclusion is reached by filtering
$\mathcal{O}_{E_j}\otimes f^*(L) ^{\otimes r}$ by line bundles on
the underlying reduced manifold. More precisely, suppose to fix
ideas that $E_j$ has generic multiplicity two, and write $E_j=2F_j$,
where $F_j$ is codimension one $G$-invariant irreducible complete
submanifold o $M$. We have an equivariant short exact sequence
$$0\longrightarrow \mathcal{O}_{F_j}\big (L^{\otimes k}(-F_j)\big )\longrightarrow
\mathcal{O}_{E_j}\otimes L^{\otimes k}\longrightarrow
\mathcal{O}_{E_j}\big ( L^{\otimes k}\big )\longrightarrow 0.$$ By
Theorem \ref{thm:pull-back} and Lemma
\ref{lem:pull-back-moment-map},
\begin{center}
$H^0\left (F_j,\mathcal{O}_{F_j}\big (L^{\otimes k}(-F_j)\big
)\right )_\mu=0$ and $H^0\left (F_j,\mathcal{O}_{F_J}\big
(L^{\otimes k}\big )\right )_\mu=0$\end{center} for every $r\gg 0$.
The general case is similar.

\bigskip

We shall need the following auxiliary result:

\begin{prop}
\label{prop:invariant-hypersurfaces-finite-vol} Let $M$ be a complex
projective $G$-manifold. Suppose that $L\in \mathrm{Pic}^G(M)$ is
ample and has non-empty stable locus: $M^\mathrm{s}(L)\neq
\emptyset$. Let $e=e_G(L)$ be the $G$-exponent of $L$ (Definition
\ref{defn:G-exponent}). Suppose $0\neq \sigma \in H^0(M,L^{\otimes
ke})^G$ is general, for some $k\gg 0$. Let
$Z=\mathrm{div}(\sigma)\subseteq M$ be the corresponding divisor.
Then $$\mathrm{vol}_\mu\left (Z,\left .A\right |_Z\right )<+\infty$$
for every $A\in \mathrm{Pic}^G(M)$ which is GIT-equivalent to $L$
and for every maximal weight $\mu$.\end{prop}

\begin{rem} Recall that ample $A,L\in \mathrm{Pic}^G(M)$ are said to be
GIT-equivalent if $M^\mathrm{u}(A)=M^{\mathrm{u}}(L)$,
$M^{\mathrm{s}}(A)=M^{\mathrm{s}}(L)$ \cite{dh}. However, the
statement of the Proposition still holds if we only require
$M^\mathrm{u}(A)\supseteq M^\mathrm{u}(L)$.
\end{rem}

\textit{Proof of Proposition
\ref{prop:invariant-hypersurfaces-finite-vol}.} We may decompose the
$G$-invariantly effective divisor $Z$ as
$Z=Z_\mathrm{u}+Z_\mathrm{s}$, where $Z_\mathrm{u},Z_{\mathrm{s}}$
are $G$-invariantly effective divisors, $Z_\mathrm{u}$ is supported
on the unstable locus $M^\mathrm{u}(L)$, and no irreducible
component of $Z_\mathrm{s}$ is supported on $M^\mathrm{u}(L)$. Since
we may assume that the base locus of $H^0(M,L^{\otimes ke})^G$ is
$M^\mathrm{u}(L)$, we may also suppose without loss that every
irreducible component of $Z_\mathrm{s}$ intersects the stable locus
$M^\mathrm{s}(L)$.

Since $Z,Z_\mathrm{u}$ and $Z_\mathrm{s}$ are all $G$-invariant
subschemes of $M$, there is a $G$-equivariant injection $
\mathcal{O}_Z\hookrightarrow \mathcal{O}_{Z_\mathrm{u}}\oplus
\mathcal{O}_{Z_\mathrm{s}}$. Tensoring by $A^{\otimes k}$, and
passing to global sections, for every maximal weight $\mu$ this
implies
$$
h^0_\mu(Z,A^{\otimes k}\otimes \mathcal{O}_{Z}) \le
h_\mu^0(Z_{\mathrm{u}},A^{\otimes k}\otimes
\mathcal{O}_{Z_{\mathrm{u}}})+ h_\mu^0(Z_{\mathrm{s}},A^{\otimes
k}\otimes \mathcal{O}_{Z_{\mathrm{s}}})$$ (where $h^0_\mu=:\dim
H^0_\mu$). Thus,
\begin{equation}\label{eqn:stable-unstable-bound}
\mathrm{vol}_\mu\left (Z,\left .A\right |_Z\right )\le
\mathrm{vol}_\mu\left (Z_\mathrm{u},\left .A\right
|_{Z_\mathrm{u}}\right )+\mathrm{vol}_\mu\left (Z_\mathrm{s},\left
.A\right |_{Z_\mathrm{s}}\right ).\end{equation}

Let us decompose $Z_\mathrm{s}$ as the sum of its irreducible
components, $Z_\mathrm{s}=\sum _jZ_{\mathrm{s}j}$. By the generality
of $\sigma$ and Bertini's Theorem, every $Z_{\mathrm{s}j}$ is
reduced and non-singular away from $M^{\mathrm{u}}(L)$. Arguing as
above, we have
\begin{equation}\label{eqn:irreduc-semist}
h_\mu^0(Z_{\mathrm{s}},A^{\otimes k}\otimes
\mathcal{O}_{Z_{\mathrm{s}}})\le \sum
_jh_\mu^0(Z_{\mathrm{s}j},A^{\otimes k}\otimes
\mathcal{O}_{Z_{\mathrm{s}j}}),\end{equation} whence
$\mathrm{vol}_\mu (Z_\mathrm{s},\left .A\right |_{Z_\mathrm{s}}) \le
\sum _j \mathrm{vol}_\mu(Z_{\mathrm{s}j},\left .A\right
|_{Z_{\mathrm{s}j}})$.

\begin{lem}\label{lem:semist-part-finite-vol-on}
$\mathrm{vol}_\mu(Z_{\mathrm{s}j},\left .A\right
|_{Z_{\mathrm{s}j}})<+\infty$ for every $j$ and $\mu$.\end{lem}

\textit{Proof of Lemma \ref{lem:semist-part-finite-vol-on}.} Let us
assume, to begin with, that the $Z_{\mathrm{s}j}$'s are
non-singular. Since $M^\mathrm{s}(\left .L\right
|_{Z_{\mathrm{s}j}})=Z_{\mathrm{s}j}\cap M^{\mathrm{s}}(L)\neq
\emptyset$ by construction, the statement follows from Corollary
\ref{cor:stable-is-enough}.

In general, let us choose for every $j$ a $G$-equivariant resolution
of singularities, $f_j:\tilde Z_{\mathrm{s}j}\rightarrow
Z_{\mathrm{s}j}$, which may be assumed to be an isomorphism away
from $Z_{\mathrm{s}j}\cap M^\mathrm{u}(L)$ \cite{ev}, \cite{eh}. We
may find a $G$-invariant effective exceptional divisor $E_j\subseteq
\tilde Z_{\mathrm{s}j}$ such that $B_j=:f_j^*(L)^{\otimes k}(-E_j)$
is an ample $G$-divisor on $\tilde{Z}_{\mathrm{s}j}$ for every $k\gg
0$. Clearly, $f_j(E_j)\subseteq Z_{\mathrm{s}j}\cap
M^\mathrm{u}(L)$. We have:

\begin{claim}\label{claim:non-empty-stable-locus-upstairs}
$M^{\mathrm{s}}(B_j)\neq \emptyset$ for all $k\gg 0$.\end{claim}

\textit{Proof of Claim \ref{claim:non-empty-stable-locus-upstairs}.}
Since $L$ is ample, section restriction $H^0(M,L^{\otimes
k})\rightarrow H^0(Z_{\mathrm{s}j},L^{\otimes k}\otimes
\mathcal{O}_{Z_{\mathrm{s}j}})$ is a $G$-equivariant surjective
linear map, for every $k\gg 0$. On the other hand, for $k$ large and
divisible, by definition we have:
$$
\mathrm{Bs}\left (\left |H^0(M,L^{\otimes k})^G\right |\right
)=M^\mathrm{u}(L).$$ Since furthermore by Nagata's Theorem the
algebra
$$R(L)^G=:\bigoplus _{k\ge 0}H^0(M,L^{\otimes k})^G$$
is finitely generated, the order of vanishing of the general $\tau
\in H^0(M,L^{\otimes k})^G$ along $M^\mathrm{u}(L)$ grows to
infinity as $k\rightarrow +\infty$; by the above, the same then
holds of the order of vanishing of the general section in
$H^0(Z_{\mathrm{s}j},L^{\otimes k}\otimes
\mathcal{O}_{Z_{\mathrm{s}j}})^G$. On the upshot, we have:

\begin{description}
  \item[i):] every invariant section of $f_j^*(L) ^{\otimes k}$ coming from
  the linear map
  \begin{equation}\label{eqn:composed-equivariant-map}H^0(M,L^{\otimes k})^G\twoheadrightarrow H^0(Z_{\mathrm{s}j},L^{\otimes k}\otimes
\mathcal{O}_{Z_{\mathrm{s}j}})^G\hookrightarrow H^0(\tilde
Z_{\mathrm{s}j},f_j^*(L)^{\otimes k})^G\end{equation} actually lives
in the image of the injection
\begin{equation*}
H^0(\tilde Z_{\mathrm{s}j},f_j^*(L)^{\otimes
k}(-E_j))^G\hookrightarrow H^0(\tilde
Z_{\mathrm{s}j},f_j^*(L)^{\otimes k})^G;\end{equation*}
  \item[ii):] when viewed as an invariant global section of $B_j=f^*(L)^{\otimes
k}(-E_j)$, every section coming from
(\ref{eqn:composed-equivariant-map}) vanishes along $E_j$.

\end{description}

Suppose now $x\in f_j^{-1}(M^{\mathrm{s}}(L))\subseteq \tilde
Z_{\mathrm{s}j}$; in particular, $x\not\in E_j$. Since the
stabilizer of $f(x)$ is discrete, so is the stabilizer of $x$. For
$k$ large and divisible, furthermore, there exists $s\in
H^0(M,L^{\otimes k})^G$ such that $s(x)\neq 0$ and the $\tilde
G$-orbit is closed in $M_s=\{p\in M:s(p)\neq 0\}$. When we pull-back
$s$ to a section of $f_j^*(L) ^{\otimes k}(-E_j)$ using i), the same
then holds in view of ii) and the fact that $f_j$ is an isomorphism
away from $M^{\mathrm{u}}(L)\cap Z_{\mathrm{s}j}$.

This completes the proof of Claim
\ref{claim:non-empty-stable-locus-upstairs}.

\bigskip

Lemma \ref{lem:semist-part-finite-vol-on} now follows in view of
Corollary \ref{cor:stable-is-enough}, since
$$\mathrm{vol}_\mu ( Z_{\mathrm{s}j},A\otimes
\mathcal{O}_{Z_{\mathrm{s}j}} )\le \mathrm{vol}_\mu (\tilde
Z_{\mathrm{s}j},f_j^*(A)).$$

\bigskip

We conclude that $\mathrm{vol}_\mu ( Z_{\mathrm{s}},A\otimes
\mathcal{O}_{Z_{\mathrm{s}}} )\le \sum _j\mathrm{vol}_\mu (
Z_{\mathrm{s}j},A\otimes \mathcal{O}_{Z_{\mathrm{s}j}} )<+\infty$.

The following is where the hypothesis that $A$ and $L$ be
GIT-equivalent is being used:

\begin{lem}\label{lem:unstable-case-finite} $\mathrm{vol}_\mu( Z_\mathrm{u},A\otimes
\mathcal{O}_{Z_\mathrm{u}})=0$ for every maximal weight
$\mu$.\end{lem}

\textit{Proof of Lemma \ref{lem:unstable-case-finite}.} It suffices
to prove that $\mathrm{vol}_\mu( Z_{\mathrm{u}j},A\otimes
\mathcal{O}_{Z_{\mathrm{u}j}})=0$ for every $j$. This follows
immediately from Proposition \ref{prop:unst-and-vanish} when
$Z_{\mathrm{u}j}$ is reduced. In the general case, we filter
$A\otimes \mathcal{O}_{Z_{\mathrm{u}j}}$ by line bundles on
$(Z_{\mathrm{u}j}) _{\mathrm{red}}$, and argue as in the final part
of the proof of Proposition \ref{prop:unst-and-vanish}.

\bigskip

Since $\mathrm{vol}_\mu( Z,A\otimes \mathcal{O}_{Z})\le
\mathrm{vol}_\mu( Z_\mathrm{s},A\otimes
\mathcal{O}_{Z_\mathrm{s}})+\mathrm{vol}_\mu( Z_\mathrm{u},A\otimes
\mathcal{O}_{Z_\mathrm{u}})$, the proof of Proposition
\ref{prop:invariant-hypersurfaces-finite-vol} is now complete.

\bigskip

The same argument also proves the following:

\begin{prop}
\label{prop:invariant-hypersurfaces-finite-vol-codimension-two} Let
$M$ be a complex projective $G$-manifold. Suppose that $L\in
\mathrm{Pic}^G(M)$ is ample and such that $M^{\mathrm{s}}(L)\neq
\emptyset$ and $\mathrm{codim}\big( M^{\mathrm{u}}(L)\big)\ge 2$.
Let $e=e_L$ be the period of $\chi _{K,L}$, as in Corollary
\ref{cor:invariant-part-asympt}. Suppose $0\neq \sigma \in
H^0(M,L^{\otimes ke})^G$ is general, for some $k\gg 0$. Let
$Z=\mathrm{div}(\sigma)\subseteq M$ be the corresponding divisor.
Then $$\mathrm{vol}_\mu\left (Z,\left .R\right |_Z\right )<+\infty$$
for \textsf{every} $R\in \mathrm{Pic}^G(M)$ and every maximal weight
$\mu$.\end{prop}

\section{Homogeneity Properties}
In this section we shall study the relation between the equivariant
volumes of $L$ and $L^{\otimes p}$, where $L\in \mathrm{Pic}^G(M)$
and $p\in \mathbb{N}$. Homogeneity does not hold in the most
simple-minded sense, as shown by Example \ref{exmp:not-homog}.

The main result here is Proposition \ref{prop:key-homog-prime}.
First however we need to establish some technical preliminaries.

In the following analysis, we may assume without loss of generality
that $M^{\mathrm{ss}}(L)\neq \emptyset$, so that
$\mathbb{N}_G(L)\neq \{0\}$, for otherwise $ \mathrm{vol}_\mu(L)=0$
for every $\mu$ in view of Proposition \ref{prop:unst-and-vanish}.
This assumption will be implicit throughout.

\begin{defn}\label{defn:volume-mod-f}
Suppose that $M$ is a complex projective $G$-variety and that $L\in
\mathrm{Pic}^G(M)$. Let $\mu$ be a maximal weight of $G$. For every
integer $f$ with $0\le f<e_G(L)$ let us define
$$
\upsilon _\mu(L,f)=:\limsup _{m\rightarrow
+\infty}\frac{h^0_\mu(M,L^{\otimes (f+me_G(L))})}{(f+me_G(L))
^{\mathrm{n-g}}/( \mathrm{n-g})!}.$$
\end{defn}

\begin{rem} $\upsilon _\mu(L,f)$ may be defined in the same manner
for any $f\in \mathbb{Z}$, and with this interpretation it only
depends on the class of $f$ mod $e_G(L)$. At places we shall refer
to this alternative interpretation, but this should cause no
confusion. \end{rem}

It is clear from the definition that $ \upsilon _\mu(L,f)\le
\mathrm{vol}_\mu(L)$ for every $f\in \{0,1,\ldots, e_G(L)-1\}$. On
the other hand, suppose that $k_i\uparrow +\infty$ is a sequence of
integers such that

$$\lim _{i\rightarrow +\infty}\frac{( \mathrm{n-g})!}{k_i^{\mathrm{n-g}}}\,h^0_\mu(M,L^{\otimes k_i})
= \mathrm{vol}_\mu(L).$$ Perhaps after passing to a subsequence, we
may assume without loss of generality that the $k_i$'s are constant
mod $e_G(L)$. Thus, if $f_0\in \{0,1,\ldots, e_G(L)-1\}$ and
$k_i\equiv f_0$ (mod. $e_G(L)$) for every $i$, then clearly $
\upsilon _\mu(L,f_0)= \mathrm{vol}_\mu(L)$. Thus,

\begin{lem} \label{lem:max-is-volume}
In the hypothesis of Definition \ref{defn:volume-mod-f},
$$ \mathrm{vol}_\mu(L)=\max \{\upsilon _\mu(L,f)\,:\,0\le
f<e_G(L)\}.$$
\end{lem}

The key to establishing the homogeneity properties of $
\mathrm{vol}_\mu$ is the following:

\begin{lem}\label{lem:f-vol-is-homog}
In the hypothesis of Definition \ref{defn:volume-mod-f}, for every
$f\in \{0,1,\ldots, e_G(L)-1\}$ and any $p\in \mathbb{N}$, we have
$$
\upsilon _\mu(L,f)=( \mathrm{n-g})!\,\limsup _{k\rightarrow
+\infty}\frac{h^0_\mu(M,L^{\otimes
(f+kp\,e_G(L))})}{\big(f+kp\,e_G(L)\big)
^{\mathrm{n-g}}}.$$\end{lem}

\textit{Proof of Lemma \ref{lem:f-vol-is-homog}.} For any integer
$r\ge 0$, let us define
$$
u_r=:( \mathrm{n}-\mathrm{g})!\,\limsup _{k\rightarrow
+\infty}\frac{h^0_\mu(M,L^{\otimes
\big(f+(kp+r)\,e_G(L))})}{\big(f+(kp+r)\,e_G(L)\big)
^{\mathrm{n-g}}}.$$ Clearly, $u_r=u_{r+p}$ for every $r$, and for
every $r_0\in \mathbb{N}$ one has $\upsilon _\mu(L,f)=\max
\{u_{r_0+1},\ldots,u_{r_0+p}\}$. We are reduced to proving the
following:

\begin{claim}\label{claim:u-r-is-constant}
There exists $r_0\in \mathrm{N}$ such that $u_r=u_0$ for all $r\in
\{r_0+1,\ldots,r_0+p\}$.\end{claim}

\textit{Proof of Claim \ref{claim:u-r-is-constant}.} Let us choose
$r_0\in \mathbb{N}$ such that for every $r\ge r_0$ we have
$H^0(M,L^{\otimes re_G(L)})^G\neq \{0\}$, for every $r\ge r_0$. We
shall argue that the statement holds for $r_0$. To see this, let us
choose $q\in \mathbb{N}$ such that $qp-(r_0+p)>r_0$. For every $r\in
\{r_0+1,\ldots,r_0+p\}$, we have $H^0(M,L^{\otimes re_G(L)})^G\neq
\{0\}$ and $H^0(M,L^{\otimes (qp-r)e_G(L)})^G\neq \{0\}$. By the
usual argument, this implies
$$
h^0_\mu\big (M,L^{\otimes \big(f+kp\,e_G(L)\big)}\big )\le
h^0_\mu\big (M,L^{\otimes \big(f+(kp+r)\,e_G(L)\big)}\big )\le
h^0_\mu\big (M,L^{\otimes \big(f+(k+q)p\,e_G(L)\big)}\big ).$$
Finally, dividing by $\big (f+kp\,e_G(L)\big ) ^{\mathrm{n-g}}/(
\mathrm{n-g})!$ and passing to $\limsup$ for $k\rightarrow +\infty$,
we are done.

\begin{cor}\label{cor:homog-of-upsilon-f}
In the hypothesis of Definition \ref{defn:volume-mod-f}, for every
$f\in \{0,1,\ldots, e_G(L)-1\}$ and any $p\in \mathbb{N}$, we have
$$
\upsilon _\mu(L^{\otimes p},f)=p^{\mathrm{n-g}}\,\upsilon _\mu
(L,pf).$$
\end{cor}

\textit{Proof.} By definition,
\begin{eqnarray}\label{eqn:homog-of-upsilon-f-1}
\upsilon _\mu(L^{\otimes p},f)&=&( \mathrm{n-g})!\,\limsup
_{m\rightarrow +\infty}\frac{h^0_\mu(M,L^{\otimes
p(f+me_G(L^{\otimes p}))})}{(f+me_G(L^{\otimes p}))
^{\mathrm{n-g}}}\nonumber
\\
&=&( \mathrm{n-g})!\,p^{\mathrm{n-g}}\,\limsup _{m\rightarrow
+\infty}\frac{h^0_\mu(M,L^{\otimes (pf+mp\,e_G(L^{\otimes
p}))})}{(pf+mp\,e_G(L^{\otimes p})) ^{\mathrm{n-g}}}.
\end{eqnarray}

By Lemma \ref{lem:e-G-L-p}, we have $e_G(L^{\otimes
p})=e_G(L)/\mathrm{gcd}(e_G(L),p)$. Let us set
$q=:p/\mathrm{gcd}(e_G(L),p)$, so that $p\,e_G(L^{\otimes
p})=q\,e_G(L)$. Inserting this in (\ref{eqn:homog-of-upsilon-f-1})
and applying Lemma \ref{lem:f-vol-is-homog} we obtain

\begin{eqnarray*}\label{eqn:homog-of-upsilon-f-2}
\upsilon _\mu(L^{\otimes p},f)&=&p^{\mathrm{n-g}}\,(
\mathrm{n-g})!\,\limsup _{m\rightarrow
+\infty}\frac{h^0_\mu(M,L^{\otimes
(pf+mq\,e_G(L))})}{(pf+mq\,e_G(L)) ^{\mathrm{n-g}}}\\
&=&p^{\mathrm{n-g}}\,\upsilon _\mu(L,p f).
\end{eqnarray*}

\bigskip

Now we can prove

\begin{prop}\label{prop:key-homog-prime}
Suppose that $L\in \mathrm{Pic}^G(M)$ and that $p\in \mathbb{N}$ is
prime with $e_G(L)$. Then $ \mathrm{vol}_\mu(L^{\otimes
p})=p^{\mathrm{n-g}}\,\mathrm{vol}_\mu(L)$.\end{prop}

\textit{Proof of Proposition \ref{prop:key-homog-prime}.} It is
immediate from the definition of equivariant volume that
$\mathrm{vol}_\mu (L^{\otimes p})\le
p^{\mathrm{n-g}}\,\mathrm{vol}_\mu (L)$.

To verify the reverse inequality, let us choose $f\in \mathbb{N}$
such that $\mathrm{vol}_\mu (L)=\upsilon _\mu(L,f)$ (Lemma
\ref{lem:max-is-volume}). Since $\upsilon _\mu(L,f)$ only depends on
the congruence class of $f$ modulo $e_G(L)$, and $p$ is prime with
$e_G(L)$ by assumption, we may assume $f=p\,\ell$ for some $\ell\in
\mathbb{N}$. Applying Corollary \ref{cor:homog-of-upsilon-f} we
conclude:

\begin{equation*}
\mathrm{vol}_\mu(L)=
\upsilon _\mu(L,p\,\ell)=
\frac{  \upsilon _\mu(L^{\otimes p},\ell)}{p^{\mathrm{n-g}}}  \le
\frac{ \mathrm{vol}_\mu(L^{\otimes p}) }{ p^{\mathrm{n-g}} }.
\end{equation*}

\begin{cor}\label{cor:general-homog}
For any $q\in \mathbb{N}$ we have
$$
\mathrm{vol}_\mu (L^{\otimes q})=\left (
\frac{q}{\mathrm{gcd}(q,e_G(L))}\right ) ^{\mathrm{n-g}}\,
\mathrm{vol}_\mu (L^{\otimes \mathrm{gcd}(q,e_G(L))}).$$
\end{cor}

\textit{Proof.} By Lemma \ref{lem:e-G-L-p}, we have
\begin{equation*}
  e_G(L^{\otimes \mathrm{gcd}(q,e_G(L))})=\frac{e_G(L)}{\mathrm{gcd}\big (e_G(L),
  \mathrm{gcd}(q,e_G(L))\big )}=\frac{e_G(L)}{
  \mathrm{gcd}(q,e_G(L))},
\end{equation*}
which is relatively prime with $q/\mathrm{gcd}(q,e_G(L))$.
Therefore, the statement follows from Proposition
\ref{prop:key-homog-prime} given that
\begin{equation}
\mathrm{vol}_\mu(L^{\otimes q})=\mathrm{vol}_\mu \left (
\big(L^{\otimes \mathrm{gcd}(q,e_G(L))}\big ) ^{\otimes
\frac{q}{\mathrm{gcd}(q,e_G(L))}}\right ).
\end{equation}

\begin{cor}\label{cor:special-case}
For any $p\in \mathbb{N}$, we have
$$ \mathrm{vol}_\mu (L^{\otimes pe_G(L)})=p^{\mathrm{n-g}}\,
\mathrm{vol}_\mu (L^{\otimes e_G(L)}).$$
\end{cor}

It is in order to single out the case of the trivial representation,
which is always homogeneous:

\begin{cor}
\label{cor:general-homog-trivial-rep} For any $L\in
\mathrm{Pic}^G(M)$ and any $q\in \mathbb{N}$ we have:
$$
\mathrm{vol}_0(L^{\otimes q})=q^{\mathrm{n-g}}\,\mathrm{vol}_0(L).
$$
\end{cor}
\textit{Proof.} If $p\in \mathbb{N}$ divides $e_G(L)$, then
\begin{eqnarray*}
\mathrm{vol}_0(L^{\otimes p})&=&( \mathrm{n-g})!\,\limsup
_{m\rightarrow +\infty}\frac{1}{m^{\mathrm{n-g}}}\,h^0(M,L^{\otimes
mp})^G\\
&=&p^{\mathrm{n-g}}\,( \mathrm{n-g})!\,\limsup _{m\rightarrow
+\infty}\frac{1}{(mp)^{\mathrm{n-g}}}\,h^0(M,L^{\otimes mp})^G\\
&=&p^{\mathrm{n-g}}\,( \mathrm{n-g})!\,\limsup _{m\rightarrow
+\infty}\frac{1}{(me_G(L))^{\mathrm{n-g}}}\,h^0(M,L^{\otimes me_G(L)})^G\\
&=&p^{\mathrm{n-g}}\,\mathrm{vol}_0(L).
\end{eqnarray*}
Here the third equality holds because $h^0(M,L^{\otimes mp})^G=0$
when $e_G(L) \nmid pm$, and the last one holds by definition of
$e_G(L)$. For any $q\in \mathbb{N}$, we then have by Corollary
\ref{cor:general-homog}:
\begin{eqnarray*}
\mathrm{vol}_0 (L^{\otimes q})&=&\left (
\frac{q}{\mathrm{gcd}(q,e_G(L))}\right ) ^{\mathrm{n-g}}\,
\mathrm{vol}_0 (L^{\otimes \mathrm{gcd}(q,e_G(L))})\\
&=&\left ( \frac{q}{\mathrm{gcd}(q,e_G(L))}\right )
^{\mathrm{n-g}}\,\mathrm{gcd}(q,e_G(L))^{\mathrm{n-g}}\,
\mathrm{vol}_0 (L)\\
&=&q^{\mathrm{\mathrm{n-g}}}\,\mathrm{vol}_0 (L),
\end{eqnarray*}
where the second equality holds because $\mathrm{gcd}(q,e_G(L))$
divides $e_G(L)$.

\section{Equivariant Kodaira Lemma}

\begin{prop}
\label{prop:kodaira-prep} Let $M$ be a complex projective
$G$-manifold, and suppose that there exists an ample $B\in
\mathrm{Pic}^G(M)$ with non-empty stable locus,
$M^{\mathrm{s}}(B)\neq \emptyset$, and such that $\mathrm{codim}\big
(M^\mathrm{u}(B)\big )\ge 2$. Then for every maximal weight $\mu$
and every $L,F\in \mathrm{Pic}^G(M)$ such that
$\mathrm{vol}_\mu(L)>0$ and $H^0(M,F)^G\neq \{0\}$, one has
$$\limsup _{k\rightarrow +\infty, k\in \mathbb{N}_\mu(L)}
\left (\frac{1}{k^{\mathrm{n-g}}}\dim H^0(M,L^{\otimes k}\otimes
F^{-1})_\mu\right )>0.$$ In particular, $ H^0(M,L^{\otimes k}\otimes
F^{-1})_\mu\neq \{0\}$ for arbitrarily large $k\in
\mathbb{N}_\mu(L)$.
\end{prop}

\textit{Proof of Proposition \ref{prop:kodaira-prep}.} To begin
with, perhaps after passing to a Kirwan resolution of the pair
$(M,B)$, we may as well assume that $B$ is ample and regular. Given
this, let us first prove the following:

\begin{lem}\label{lem:kod-auxiliary}
Let $e=e_G(B)$ be the $G$-exponent of $B$. Then
$$H^0(M,B^{\otimes re_G(B)}\otimes F^{-1})^G\neq \{0\}$$
if $r\gg 0$.
\end{lem}

\textit{Proof of Lemma \ref{lem:kod-auxiliary}.} Let $(K)$ denote
the principal orbit type of the $G$-action on $M$, and let $\chi
_{K,F}:K\rightarrow S^1$ be the character associated to $F\in
\mathrm{Pic}^G(M)$, as in Remarks
\ref{rem:condtion-independ-of-p-chi} and \ref{rem:principal-type}.
The hypothesis $H^0(M,F)^G\neq \{0\}$ implies $\chi _{K,F}=1$.
Therefore, the character associated to $B^{\otimes re_G(B)}\otimes
F^{-1}\in \mathrm{Pic}^G(M)$ is $(\chi_{K,B})^{re_G(B)}\cdot \chi
_{K,F}^{-1}=(\chi_{K,B})^{re_G(B)}=1$.

On the other hand, for $r\gg 0$, $ B^{\otimes re_G(B)}\otimes
F^{-1}\in \mathrm{Pic}^G(M)$ is ample and regular, as its
equivalence class in $ \mathrm{NS}^G(M)_\mathbb{R}$ lies in the
interior of the same chamber as the class of $B$. Therefore, arguing
as in the proof of Theorem \ref{thm:main-regular-case}, the
dimension of $H^0(M,B^{\otimes re_G(B)}\otimes F^{-1})^G$ may be
computed using the Riemann-Roch formula for multiplicities due to
Meinrenken.

Now the moment map associated to $B^{\otimes re}\otimes F^{-1}$ is
$$\Phi
_{B^{\otimes re}\otimes F^{-1}}=:re\Phi_B-\Phi_F:M\rightarrow
\frak{g}^*,$$ where $\Phi _F$ is a moment map for $F$. Consequently,
$\Phi _{B^{\otimes re}\otimes F^{-1}}^{-1}(0)\rightarrow \Phi
_B^{-1}(0)$, in the following sense: For $\epsilon \in \mathbb{R}$
sufficiently small, let us set $\Phi _\epsilon = :\Phi _B-\epsilon
\Phi _F$. Since $0\in \frak{g}^*$ is a regular value of $\Phi
_0=\Phi _B$, and $\Phi _0^{-1}(0)\neq \emptyset$, the same holds of
$\Phi _\epsilon$ for sufficiently small $\epsilon$. Thus, the loci
$\Phi _\epsilon ^{-1}(0)\subseteq M$ form a family of compact
connected g-codimensional submanifolds of $M$. The statement of the
Lemma now follows from the same asymptotic computations as in the
proof of Theorem \ref{thm:main-regular-case} and Corollary
\ref{cor:main-regular-case-but-general}.

\bigskip

Given the Lemma, let us choose $0\neq \tau \in H^0(M,B^{\otimes
re}\otimes F^{-1})^G$ for some $r\gg 0$. Tensor product by $\tau$
determines injective linear maps
$$
H^0(M,L^{\otimes k}\otimes B^{-re})_\mu\hookrightarrow
H^0(M,L^{\otimes k}\otimes F^{-1})_\mu$$ for every $k$ and $\mu$. We
are thus reduced to proving the following:

\begin{claim}\label{claim:asymptotic-non-vanishing-at-mu-F}
For every maximal weight $\mu$ and every integer $r\gg 0$, there
exists a sequence $k_i\in \mathbb{N}_\mu(L)$, $k_i\uparrow +\infty$
such that
\begin{equation}\label{eqn:non-vanishing}
\dim H^0(M,L^{\otimes k_i}\otimes
B^{-re})_\mu=O(k_i^{\mathrm{n-g}}).
\end{equation}
\end{claim}

\textit{Proof of Claim
\ref{claim:asymptotic-non-vanishing-at-mu-F}.} Since
$+\infty>\mathrm{vol}_\mu(L)>0$, we can find a sequence $k_i\uparrow
+\infty$, $k_i\in \mathbb{N}_\mu(L)$, such that
\begin{equation}\label{eqn:ses-equiv-sigma}
\dim H^0(M,L^{\otimes k_i})_\mu=O(k_i^{\mathrm{n-g}}).
\end{equation}
For $r\gg 0$ and divisible, let us choose a general $\sigma \in
H^0(M,B^{re})^G$, and let $Z_\sigma=\mathrm{div}(\sigma)$ be its
zero divisor. For every $i$ we have a $G$-equivariant short exact
sequence
\begin{equation*}
0\longrightarrow L^{\otimes k_i}\otimes B^{-re}\longrightarrow
L^{\otimes k_i}\longrightarrow \mathcal{O}_{Z_\sigma}(L^{\otimes
k_i})\longrightarrow 0.
\end{equation*}
Passing to the $\mu$-th equivariant summand in cohomology, we obtain
the other
\begin{equation}\label{eqn:ses-equiv-sigma-coho}
0\rightarrow H^0(M,L^{\otimes k_i}\otimes B^{-re})_\mu\rightarrow
H^0(M,L^{\otimes k_i})_\mu\rightarrow
H^0(Z_\sigma,\mathcal{O}_{Z_\sigma}(L^{\otimes k_i}))_\mu.
\end{equation}
In view of Proposition
\ref{prop:invariant-hypersurfaces-finite-vol-codimension-two}, we
have $\dim H^0(Z_\sigma,\mathcal{O}_{Z_\sigma}(L^{\otimes
k_i}))_\mu=o(k_i^{\mathrm{n-g}})$. The statement of Claim
\ref{claim:asymptotic-non-vanishing-at-mu-F} follows from
(\ref{eqn:ses-equiv-sigma-coho}) in view of
(\ref{eqn:ses-equiv-sigma}).

\bigskip

\begin{cor}\label{cor:vol-mu-pos-equiv-charat}
Let $M$ be a complex projective $G$-manifold, and suppose that there
exists an ample $B\in \mathrm{Pic}^G(M)$ such that $M^{s}(B)\neq
\emptyset$, $\mathrm{codim}\big (M^\mathrm{u}(B)\big)\geq 2$. Then
the following conditions on $L\in \mathrm{Pic}^G(M)$ are equivalent:
\begin{description}
  \item[i):] $ \mathrm{vol}_\mu(L)>0$;
  \item[ii):] for every $F\in \mathrm{Pic}^G(M)$ with $
  \mathrm{vol}_0(F)>0$, there exists a sequence $k_i\in
  \mathbb{N}_\mu(L)$, $k_i\uparrow +\infty$, such that $$L^{\otimes k_i}=F^{\otimes
  e_G(F)}\otimes A_i,$$ where
  $\dim H^0(M,A_i)_\mu=O(k_i^{\mathrm{n-g}})$.
\end{description}
\end{cor}

\textit{Proof of Corollary \ref{cor:vol-mu-pos-equiv-charat}.} That
i) implies ii) is the content of Proposition
\ref{prop:kodaira-prep}. The reverse implication is obvious.

\begin{cor}\label{cor_G-big-implies-big}
In the hypothesis of Corollary \ref{cor:vol-mu-pos-equiv-charat}, if
$L\in \mathrm{Pic}^G(M)$ satisfies $ \mathrm{vol}_0(L)>0$, then the
underlying line bundle $f_G(L)\in \mathrm{Pic}(M)$ is big.\end{cor}

\begin{exmp}\label{exmp:codim-is-important}
The following example shows that in general the hypothesis of
Corollary \ref{cor:vol-mu-pos-equiv-charat} may not be replaced by
the weaker one $M^\mathrm{s}(B)\neq \emptyset$: Let $S^1$ act on $
\mathbb{C}^2$ by $t\cdot (z_0,z_1)=:(tz_0,t^{-1}z_1)$, and consider
the induced action on $ \mathbb{P}^1$, with the built-in
linearization to the hyperplane bundle, $H\in \mathrm{Pic}^{S^1}(
\mathbb{P}^1)$. Let us consider the product action of $S^1\times
S^1$ on $ \mathbb{P}^1\times \mathbb{P}^1$. Let $\pi
_i:\mathbb{P}^1\times \mathbb{P}^1\rightarrow \mathbb{P}^1$ be the
projections, and define $H_i=\pi _i^*(H)\in
\mathrm{Pic}^G(\mathbb{P}^1\times \mathbb{P}^1)$. Then
$\mathrm{vol}_0(H_i)=\mathrm{vol}_0(H)>0$, but the underlying line
bundle of $H_i$ is the pull-back by $\pi _i$ of the hyperplane
bundle on $ \mathbb{P}^1$, which is not big.
\end{exmp}


\begin{cor} \label{lem:G-big-implies-eG-is-nice}
In the situation of Corollary \ref{cor:vol-mu-pos-equiv-charat}
suppose $L\in \mathrm{Pic}^G(M)$ satisfies $\mathrm{vol}_0(L)>0$.
Let $(K)$ be the principal orbit type of the action. Then $e_G(L)$
equals the period of $\chi _{K,L}$, $|\chi _{K,L}|$.
\end{cor}


\textit{Proof of Corollary \ref{lem:G-big-implies-eG-is-nice}.} The
following argument is inspired by the proof of Corollary 2.2.10 of
\cite{lazI}.

If $H^0(M,L^{\otimes k})^G\neq \{0\}$, then $\chi _{K,L^{\otimes
k}}=(\chi _{K,L}) ^k=1$. Thus, every sufficiently large multiple of
$e_G(L)$ is a multiple of $|\chi _{K,L}|$ so that $|\chi _{K,L}|$
divides $e_G(L)$.

To prove that, conversely, $e_G(L)$ divides $|\chi _{K,L}|$, let us
replace $L$ by $L^{\otimes |\chi _{K,L}|}$, so that $|\chi
_{K,L}|=1$; by Lemma \ref{lem:e-G-L-p} it suffices to prove that
$e_G(L)=1$.

Perhaps after passing to a Kirwan resolution of the pair $(M,B)$, we
may assume without loss that $B$ is ample and regular. By Lemma
\ref{lem:kod-auxiliary} (or its proof), perhaps after replacing $B$
by $B^{\otimes re_G(B)}$ we may assume that $D=:B\otimes L^{-1}$ is
$G$-invariantly effective. Since on the other hand $
\mathrm{vol}_0(L)>0$, by Corollary \ref{cor:vol-mu-pos-equiv-charat}
we may find $m\gg 0$ such that $L^{\otimes m}$ and $E=:L^{\otimes
m}\otimes B^{-1}$ are both $G$-invariantly effective (recall the
choice of the exponents $k_i$ described in
(\ref{eqn:ses-equiv-sigma})). Thus,
$$
L^{\otimes (m-1)}=L^{\otimes m}\otimes L^{-1}=L^{\otimes m}\otimes
D\otimes B^{-1}=D\otimes E$$ is also $G$-invariantly effective. It
follows that $e_G(L)=1$.



\begin{rem}\label{rem:locally-free-case}
Let us note the following obvious special case. In the situation of
Corollary \ref{cor:vol-mu-pos-equiv-charat} assume in addition that
the action of $G$ on $M$ is generically free, so that $K$ is
trivial: $K=(e)$. If $L\in \mathrm{Pic}^G(M)$ satisfies
$\mathrm{vol}_0(L)>0$, then $e_G(L)=1$.
\end{rem}

In view of Remark \ref{rem:locally-free-case} and Corollary
\ref{cor:special-case}, we obtain:


\begin{cor}\label{cor:big-and-locally-free-homogeneous}
Suppose that the action of $G$ on $M$ is generically free, and that
the hypothesis of Corollary \ref{cor:vol-mu-pos-equiv-charat} are
satisfied. If $\mathrm{vol}_0(L)>0$, then
$\mathrm{vol}_\mu(L^{\otimes
p})=p^{\mathrm{n-g}}\,\mathrm{vol}_\mu(L)$ for every highest weight
$\mu$ and every $p=1,2,\ldots$.
\end{cor}


\begin{lem} \label{lem:passage-to-divisors-bound}
Suppose that $B\in \mathrm{Pic}^G(M)$ is ample and such that
$M^{s}(B)\neq \emptyset$, $\mathrm{codim}\big
(M^\mathrm{u}(B)\big)\geq 2$. Choose $r\gg 0$ and a very general
$\sigma \in H^0(M,B^{\otimes re_G(B)})^G$. Then for every $L\in
\mathrm{Pic}^G(M)$ and $m\in \mathbb{N}$ we have
\begin{eqnarray*} h^0_\mu\left (M,(L\otimes
B^{-re_G(B)})^m\right ) \ge h^0_\mu(M,L^{\otimes
m})-m\,h_\mu^0(Z,L^{\otimes m}\otimes \mathcal{O}_Z),
\end{eqnarray*}
where $Z=\mathrm{zero}(\sigma)$ is the zero locus of $\sigma$ (here
$h^0_\mu=\dim H^0_\mu$).\end{lem}


\textit{Proof of Lemma \ref{lem:passage-to-divisors-bound}.} If
$r\gg 0$, we have $\mathrm{Bs}\left (\left |H^0(M,B^{\otimes
re_G(B)})^G\right |\right )=M^\mathrm{u}(B)$, and by assumption the
latter has codimension $\ge 2$. Therefore, for general $\sigma
_1,\ldots, \sigma _m\in H^0(M,B^{\otimes re_G(B)})^G$ the zero loci
$Z_j=\mathrm{zero}(\sigma _j)$ have no irreducible component in
common. Hence, tensor power by $\sigma _1\otimes \cdots \otimes
\sigma _m$ induces an exact sequence
$$0\rightarrow H^0\left (M,(L\otimes
B^{-re_G(B)})^m\right )_\mu\rightarrow H^0\left (M,L^{\otimes
m}\right )_\mu\rightarrow \bigoplus _{j=1}^m H^0\left
(Z_j,L^m\otimes \mathcal{O}_{Z_j} \right )_\mu.$$ Given a very
general $\sigma \in H^0(M,B^{\otimes re_G(B)})^G$, we deduce from
this
\begin{eqnarray*} h^0_\mu\left (M,(L\otimes
B^{-re_G(B)})^m\right ) &\ge& h^0_\mu(M,L^{\otimes m})-\sum
_jh_\mu^0(Z_j,L^{\otimes m}\otimes \mathcal{O}_{Z_j})\\
&= & h^0_\mu(M,L^{\otimes m})-m\cdot h_\mu^0(Z,L^{\otimes m}\otimes
\mathcal{O}_Z),
\end{eqnarray*}
where the latter equality follows from semicontinuity and the very
generality of $\sigma$.


\begin{prop}
\label{prop:limit-infinity-p-prime} Let $M$ be a complex projective
$G$-manifold, and suppose that $ B\in\mathrm{Pic}^G(M)$ is ample and
such that $M^{s}(B)\neq \emptyset$, $\mathrm{codim}\big
(M^\mathrm{u}(B)\big)\geq 2$. Then
$$
\limsup _{p\rightarrow +\infty}\frac{\mathrm{vol}_\mu(L^{\otimes p
}\otimes B^{\otimes (-r
e_G(B))})}{p^{\mathrm{n-g}}}=\mathrm{vol}_\mu(L)$$ for every $L\in
\mathrm{Pic}^G(M)$, $r=1,2,\ldots$ and maximal weight $\mu$.
\end{prop}


\textit{Proof.} To simplify notation, let us write $C=B^{\otimes
re_G(B) }$. We may assume without loss that $r\gg 0$, so that
$\mathrm{Bs}(H^0(M,C)^G)=\mathrm{M}^{\mathrm{u}}(B)$. Since $C$ is
$G$-invariantly effective, by Remark \ref{rem:G-eff-morphism} we
have
\begin{equation}\label{eqn:vol-mu-bound-p}
 \mathrm{vol}_\mu(L^{\otimes p
}\otimes C^{-1}))\le \mathrm{vol}_\mu(L^{\otimes p})\le
p^{\mathrm{n-g}}\,\mathrm{vol}_\mu (L);\end{equation} the second
inequality in (\ref{eqn:vol-mu-bound-p}) is immediate from the
definition of volume.

On the other hand, if $Z=\mathrm{zero}(\sigma)$ for some very
general $\sigma \in H^0(M,C)^G$, by Lemma
\ref{lem:passage-to-divisors-bound} applied with $L$ replaced by
$L^{\otimes p}$ we have
\begin{equation*}
h^0_\mu\left (M,(L^{\otimes p}\otimes C^{-1})^{\otimes m}\right )
\ge h^0_\mu(M,L^{\otimes pm})-m\,h_\mu^0(Z,L^{\otimes pm}\otimes
\mathcal{O}_Z).
\end{equation*}
Dividing by $m^{\mathrm{n-g}}/( \mathrm{n-g})!$ and taking the
$\limsup$, we obtain:
\begin{equation*}
\mathrm{vol}_\mu (M,L^{\otimes p}\otimes C^{-1})\ge \mathrm{vol}_\mu
(M,L^{\otimes p})-(\mathrm{n-g})\, \mathrm{vol}_\mu (Z,L^{\otimes
p}\otimes \mathcal{O}_Z)
\end{equation*}
Now we remark that $\mathrm{vol}_\mu (Z,L^{\otimes p}\otimes
\mathcal{O}_Z)\le p^{\mathrm{n-g-1}}\,\mathrm{vol}_\mu (Z,L\otimes
\mathcal{O}_Z)$, and that the latter volume is finite by Proposition
\ref{prop:invariant-hypersurfaces-finite-vol-codimension-two}.
Therefore, dividing by $p^{\mathrm{n-g}}$ and taking the $\limsup$,
we obtain:
\begin{equation*}
\limsup _{p\rightarrow +\infty }\frac{\mathrm{vol}_\mu (M,L^{\otimes
p}\otimes C^{-1})}{p^{\mathrm{n-g}}}\ge \limsup _{p\rightarrow
+\infty }\frac{\mathrm{vol}_\mu (M,L^{\otimes
p})}{p^{\mathrm{n-g}}}.
\end{equation*}
The statement of Proposition \ref{prop:limit-infinity-p-prime} is
now a consequence of Proposition \ref{prop:key-homog-prime}.

\bigskip

If in addition $B\in \mathrm{Pic}^G(M)$ as in the hypothesis of
Proposition \ref{prop:limit-infinity-p-prime} can be chosen regular,
a stronger statement holds:

\begin{cor}\label{cor:limit-infinity-p-prime}
Suppose that there exists $ B\in\mathrm{Pic}^G(M)$ ample and regular
and such that $\mathrm{codim}\big (M^\mathrm{u}(B)\big)\geq 2$. Let
$A\in \mathrm{Pic}^G(M)$ be any ample $G$-linearized line bundle
with $M^\mathrm{s}(A)\neq \emptyset$. Suppose $r\gg 0$. Then
$$
\limsup _{p\rightarrow +\infty}\frac{\mathrm{vol}_\mu(L^{\otimes p
}\otimes A^{\otimes (-r
e_G(A))})}{p^{\mathrm{n-g}}}=\mathrm{vol}_\mu(L)$$ for every $L\in
\mathrm{Pic}^G(M)$ and maximal weight $\mu$.
\end{cor}

\textit{Proof.} One inequality is obvious, as in the proof of
Proposition \ref{prop:limit-infinity-p-prime}. In the opposite
direction, let us choose $s\gg 0$ so that the numerical equivalence
class of $H=A\otimes B^{\otimes se_G(B)}$ lies in the interior of
the same chamber as the class of $B$; in particular, $H$ is ample
and regular, and $\mathrm{codim}\big ( M^\mathrm{u}(H)\big)\ge 2$.
We may assume that $B^{\otimes nse_G(B)}$ is effective for every
integer $n\gg 0$. We have $e_G(H)=e_G(A)$. Therefore, setting
$r'=rse_G(B)$, we have:
\begin{eqnarray*}
\limsup_{p\rightarrow +\infty}\frac{\mathrm{vol}_\mu \big(L^{\otimes
p}\otimes
A^{-re_G(A)}\big)}{p^{\mathrm{n-g}}}&\ge&\limsup_{p\rightarrow
+\infty}\frac{\mathrm{vol}_\mu \big(L^{\otimes p}\otimes (A\otimes
B^{\otimes se_G(B)}) ^{-re_G(A)}\big)}{p^{\mathrm{n-g}}}\\
&=&\limsup_{p\rightarrow +\infty}\frac{\mathrm{vol}_\mu
\big(L^{\otimes p}\otimes H ^{-r'e_G(H)}\big)}{p^{\mathrm{n-g}}}\\
&=&\mathrm{vol}_\mu(L).
\end{eqnarray*}
The last equality holds by Proposition
\ref{prop:limit-infinity-p-prime}.

\begin{rem}\label{rem:lim-sup-and-lim}
In view of Proposition \ref{prop:key-homog-prime}, one can actually
strengthen the statement of Proposition
\ref{prop:limit-infinity-p-prime} (and Corollary
\ref{cor:limit-infinity-p-prime}) as follows: let $p_i\uparrow
+\infty$ be a sequence of positive integers prime with $e_G(L)$.
Then
$$
\lim _{i\rightarrow +\infty}\frac{\mathrm{vol}_\mu(L^{\otimes p_i
}\otimes B^{\otimes (-r
e_G(B))})}{p_i^{\mathrm{n-g}}}=\mathrm{vol}_\mu(L).$$
\end{rem}

The equivariant Kodaira Lemma (Corollary
\ref{cor:vol-mu-pos-equiv-charat}) implies a characterization of nef
and big line bundles with positive 0-volume.

\begin{cor}\label{cor:nef-and-big-case}
Suppose that there exists $ B\in\mathrm{Pic}^G(M)$ ample and such
that $M^{s}(B)\neq \emptyset$, $\mathrm{codim}\big
(M^\mathrm{u}(B)\big)\geq 2$. Then the following conditions are
equivalent on $L\in \mathrm{Pic}^G(M)$:
\begin{enumerate}
    \item  $L$ is nef and big, and $\mathrm{vol}_0(L)>0$;
    \item  there exists $E\in \mathrm{Pic}^G(M)$ which is
    $G$-invariantly effective and is such that $L^{\otimes
    k}\otimes E^{-1}\in \mathrm{Pic}^G(M)$ is ample and has
    non-empty stable locus, for every $k\gg 0$.
\end{enumerate}
\end{cor}

\textit{Proof.} 1. implies 2.: Perhaps after replacing $B$ by some
power $B^{\otimes re_G(B)}$, $r\gg 0$, we may assume that
$H^0(M,B)^G$ globally generates $B$ on the dense open subset
$M^{ss}(B)$. By Corollary \ref{cor:vol-mu-pos-equiv-charat}, we may
find $k_0\gg 0$ and a $G$-invariantly effective $E\in
\mathrm{Pic}^G(M)$ such that $L^{\otimes k_0}=B\otimes E$. Suppose
$k\gg 0$ and $H^0(M,L^{\otimes k})^G\neq \{0\}$, and choose $\tau
\in H^0(M,L^{\otimes k})^G$, $\tau \neq 0$. Let $U_\tau=\{p\in
M:\tau (p)\neq 0\}$ and choose $p\in U_\tau \cap M^\mathrm{s}(B)$:
in particular, $p$ has finite stabilizer $\widetilde{G}_p\subseteq
\widetilde{G}$. We can find $\sigma \in H^0(M,B)^G$ such that
$\sigma (p)\neq 0$ and the $\widetilde{G}$-orbit of $p$ is closed in
$U_\sigma =\{p:\sigma (p)\neq 0\}$. Then clearly $\tau \otimes
\sigma (p)\neq 0$ and $\widetilde{G}\cdot p$ is closed in
$U_{\tau\otimes \sigma}\subseteq U_\sigma$. Thus $L^{\otimes
k}\otimes B$ is ample and $p\in M^\mathrm{s}\Big (L^{\otimes
k}\otimes B\Big )$. Now $L^{\otimes (k+k_0)}=\Big (L^{\otimes
k}\otimes B\Big )\otimes E$.

2. implies 1.: obvious.

\section{Relation to numerical equivalence}

We shall now show that, at least under certain hypothesis, the
equivariant volumes $\mathrm{vol}_\mu(L)$ depend only on the
numerical equivalence class of $L\in \mathrm{Pic}^G(M)$.

Referring to Definition \ref{defn:homolo-trivial-G} and Remark
\ref{rem:jacobian-connected}, we have:

\begin{lem} \label{lem:homolo-equiv-zero-effective}
Suppose $B\in \mathrm{Pic}^G(M)$ is ample and $M^\mathrm{s}(B)\neq
\emptyset$. Then there exist arbitrarily large integers $r\in
\mathbb{N}$ such that
$$
H^0(M,B^{\otimes re_G(B)}\otimes P)^G\neq \{0\},$$ for every $P\in
\mathrm{Pic}^{G}(M)_0'$.
\end{lem}

\textit{Proof of Lemma \ref{lem:homolo-equiv-zero-effective}.} To
begin with, let us remark that if $f:\widetilde{M}\rightarrow M$ is
an equivariant projective morphism, then the pull-back
$f^*:\mathrm{Pic}^G(M)\rightarrow \mathrm{Pic}^G( \widetilde{M})$
satisfies $$f^*\left (\mathrm{Pic}^G(M)_0\right )\subseteq
\mathrm{Pic}^G(\widetilde{M})_0,\,\,f^*\left
(\mathrm{Pic}^G(M)_0'\right )\subseteq
\mathrm{Pic}^G(\widetilde{M})_0'.$$ Thus, without loss of generality
we may pass if necessary to a Kirwan resolution of the pair $(M,B)$
so as to assume that $B$ is ample and regular. Let us also replace
$B$ by $B^{\otimes e_G(B)}$, so as to assume without loss that
$e_G(B)=1$.

For every $P\in \mathrm{Pic}^G(M)_0'$, the following holds
(\cite{dh}, Propositions 2.3.3 and 3.1.4):
\begin{description}
  \item[i):] $M^{\mathrm{ss}}(B\otimes P)=M^\mathrm{ss}(B)\neq
\emptyset$;
  \item[ii):] $e_G(B)=e_G(B\otimes P)=1$;

  therefore,
\item[iii):] there exists $r_P$ such that  $H^0\left (M,(B\otimes
P) ^{\otimes r}\right )^G\neq \{0\}$ for all $r\ge r_P$.
\end{description}
Hence, by semicontinuity and compactness, there exists $r_0$ such
that
$$H^0\left (M,B ^{\otimes r}\otimes P ^{\otimes
r}\right )^G=H^0\left (M,(B\otimes P) ^{\otimes r}\right )^G\neq
\{0\}$$ for all $P\in \mathrm{Pic}^G(M)_0'$ and $r\ge r_0$. Now
suppose that $ H^2(M,\mathbb{Z}) _{\mathrm{tor}}\cong
\mathbb{Z}_{a_1}^{\oplus m_1}\oplus \cdots \oplus
\mathbb{Z}_{a_k}^{\oplus m_k}$ for certain integers $a_i,m_i\in
\mathbb{N}$. If $r\ge r_0$ is prime with every $a_i$, then $P\mapsto
P^{\otimes r}$ is a surjective morphism $
\mathrm{Pic}^G(M)_0'\rightarrow \mathrm{Pic}^G(M)_0'$. The statement
follows.

\begin{cor} \label{cor:homolo-equiv-ok}
Suppose that there exists an ample $B\in \mathrm{Pic}^G(M)$ such
that $M^\mathrm{s}(B)\neq \emptyset$. If $L,L'\in \mathrm{Pic}^G(M)$
are numerically equivalent and $\mu$ is a maximal weight, then
$\mathrm{vol}_\mu(L)>0$ if and only
$\mathrm{vol}_\mu(L')>0$.\end{cor}

\textit{Proof.} Perhaps after passing to a Kirwan resolution, we may
assume without loss that $B$ is ample and regular. By Lemma
\ref{lem:homolo-equiv-zero-effective}, perhaps after replacing $B$
by $B^{\otimes re_G(B)}$ for some $r\gg 0$ we may assume without
loss that $e_G(B)=1$ and that $H^0(M,B\otimes P)^G\neq \{0\}$ for
every $P\in \mathrm{Pic}^G(M)_0'$.

Suppose $\mathrm{vol}_\mu(L)>0$. By Corollary
\ref{cor:vol-mu-pos-equiv-charat}, there exists a sequence $k_i\in
\mathbb{N}_\mu(L)$ with $k_i\uparrow +\infty$ such that $L^{\otimes
k_i}=B\otimes F_i$, with $\dim H^0(M,F_i)_\mu=O(k_i^{\mathrm{n-g}})$
  as $i\rightarrow +\infty$.

  Since $L\sim _\mathrm{n} L'$, we have $L'=L\otimes P$ for some
  $P\in \mathrm{Pic}^G(M)'_0$.
    Thus,
  $(L') ^{\otimes k_i}=(B\otimes P^{\otimes k_i}) \otimes F_i$.
Now $H^0(M,B\otimes P^{\otimes k_i})^G\neq \{0\}$, and therefore
$\dim H^0(M,F_i)_\mu \le \dim H^0(M,(L') ^{\otimes k_i})_\mu$.
Consequently, $ \dim H^0(M,(L') ^{\otimes
k_i})_\mu=O(k_i^{\mathrm{n-g}})$, so that $\mathrm{vol}_\mu(L')>0$.

\begin{thm}\label{thm:vol-indep-homolo-equiv}
Let us suppose that there exists $B\in \mathrm{Pic}^G(M)$ ample and
such that $M^\mathrm{s}(B)\neq \emptyset$, $\mathrm{codim}\left
(M^\mathrm{u}(B)\right )\ge 2$. If $L,L'\in \mathrm{Pic}^G(M)$,
$L\sim _\mathrm{n} L'$ and $\mu$ is a maximal weight, then
$$ \mathrm{vol}_\mu (L)=\mathrm{vol}_\mu (L').$$
\end{thm}

\textit{Proof of Theorem \ref{thm:vol-indep-homolo-equiv}.} By
assumption, $L'=L\otimes P$ for some $P\in \mathrm{Pic}^G(M)_0'$.
Given Corollary \ref{cor:homolo-equiv-ok}, we may assume that $
\mathrm{vol}_\mu (L)>0$. In view of Lemma
\ref{lem:homolo-equiv-zero-effective}, there exists $r\in
\mathbb{N}$ such that $H^0(M,B^{\otimes re_G(B)}\otimes
P^{-p})^G\neq \{0\}$ for every $p\in \mathbb{Z}$. Recalling Remark
\ref{rem:G-eff-morphism}, we then have
$$ \mathrm{vol}_\mu \Big ((L\otimes P) ^{\otimes p}\otimes B^{-\otimes re_G(B)}
\Big )\le \mathrm{vol}_\mu \big(L^{\otimes p}\big)\le
p^{\mathrm{n-g}} \,\mathrm{vol}_\mu (L),$$ for any maximal weight
$\mu$. Proposition \ref{prop:limit-infinity-p-prime}, applied with
$L$ replaced by $L\otimes P$, now implies $ \mathrm{vol}_\mu \left
(L\otimes P \right )\le \,\mathrm{vol}_\mu (L)$. The statement then
follows by exchanging the roles of $P$ and $P^{-1}$.


\begin{rem}\label{rem:well-defined-over-Q}
Given Theorem \ref{thm:vol-indep-homolo-equiv}, for every maximal
weight $\mu$ there is a well-defined function
$\mathrm{vol}_\mu:\mathrm{Num}^G(M)\rightarrow \mathbb{R}$. In view
of Example \ref{exmp:not-homog}, this does not extend to a
well-defined function on $\mathrm{Num}^G(M)
_\mathbb{Q}=:\mathrm{Num}^G(M)\otimes \mathbb{Q}$ in a natural way,
except in the notable case $\mu =0$. In fact, momentarily adopting
additive notation, given any $P\in \mathrm{NS}^G(M) _\mathbb{Q}$,
there exists $a\in \mathbb{N}$ such that $aP\in \mathrm{NS}^G(M)$;
we may then let $
\mathrm{vol}_0(P)=:\frac{1}{a^{\mathrm{n-g}}}\mathrm{vol}_0(aP)$.
This is well-defined (that is, independent of the choice of $a$) by
Corollary \ref{cor:general-homog-trivial-rep}.
\end{rem}

\section{The $G$-big cone}

Recall that $\mathrm{C}^G(M)\subseteq \mathrm{NS}^G(M)_\mathbb{R}$
is the $G$-ample cone of \cite{dh}. The integral points in its
interior $\mathrm{Int}\left (\mathrm{C}^G(M)\right )$ are the
numerical equivalence classes of the ample $L\in \mathrm{Pic}^G(M)$
such that $M^\mathrm{s}(L)\neq \emptyset$ (the non-rational points
$\xi \in \mathrm{Int}\left (\mathrm{C}^G(M)\right )$ have a similar
interpretation in terms of K\"{a}hler classes \cite{dh}). Thus, in
the situation of Corollary \ref{cor:vol-0-pos-equiv-charat} below,
$\mathrm{Int}\left (\mathrm{C}^G(M)\right )$ is a non-empty open
cone in $\mathrm{NS}^G(M)_\mathbb{R}$.

In the special case $\mu=0$, Corollary
\ref{cor:vol-mu-pos-equiv-charat} yields the following
characterization:

\begin{cor}\label{cor:vol-0-pos-equiv-charat}
Let $M$ be a complex projective $G$-manifold, and suppose that there
exists an ample $B\in \mathrm{Pic}^G(M)$ such that
$M^\mathrm{s}(B)\neq \emptyset$, $\mathrm{codim}\left
(M^\mathrm{u}(B)\right )\ge 2$. Then the following conditions on
$L\in \mathrm{Pic}^G(M)$ are equivalent:
\begin{description}
  \item[i):] $ \mathrm{vol}_0(L)>0$;
  \item[ii):] there exist $k\in \mathbb{N}$ and $F\in \mathrm{Pic}^G(M)$
  with
  $H^0(M,F)^G\neq \{0\}$ such that $L^{\otimes
  k}=B^{e_G(B)}\otimes F$;
  \item[iii):] there exist $k\in \mathbb{N}$, $A,\,F\in \mathrm{Pic}^G(M)$
  with $A\in \mathrm{Int}\left (\mathrm{C}^G(M)\right )$,
  $H^0(M,F)^G\neq \{0\}$, such that $L^{\otimes
  k}=B^{e_G(B)}\otimes F$.
\end{description}
\end{cor}

\begin{defn} \label{defn:G-big-G-eff}
\begin{description}
    \item[i):] An element $\xi \in \mathrm{NS}^G(M)_\mathbb{R}$ will be called
\textit{$G$-big} if it has the form $\xi =\sum _{i=1}^ra_i[L_i]$,
where $a_i>0$ and $L_i\in \mathrm{Pic}^G(M)$ satisfies $
\mathrm{vol}_0(L_i)>0$ for every $i$. The \textit{$G$-big cone}
$\mathrm{Big}^G(M)\subseteq \mathrm{NS}^G(M)_\mathbb{R}$ is the
convex cone of all $G$-big elements of
$\mathrm{NS}^G(M)_\mathbb{R}$.
    \item[ii):] An element $\xi \in \mathrm{NS}^G(M)_\mathbb{R}$ will be called
\textit{$G$-effective} if it has the form $\xi =\sum
_{i=1}^ra_i[L_i]$, where $a_i>0$ and $L_i\in \mathrm{Pic}^G(M)$
satisfies $M^\mathrm{ss}(L_i)\neq \emptyset$ for every $i$. The
\textit{$G$-effective cone} $\mathrm{Eff}^G(M)\subseteq
\mathrm{NS}^G(M)_\mathbb{R}$ is the convex cone of all $G$-effective
elements of $\mathrm{NS}^G(M)_\mathbb{R}$.
\item[iii):] The \textit{$G$-pseudo-effective cone} is the closure
$\overline{\mathrm{Eff}^G(M)}\subseteq \mathrm{NS}^G(M)_\mathbb{R}$.
\end{description}
\end{defn}

By analogy with Proposition 2.2.22 of \cite{lazI}, and essentially
by the same argument, we have:

\begin{prop}\label{prop:big-R-classes-charact}
In the situation of Corollary \ref{cor:vol-0-pos-equiv-charat},
$\xi\in \mathrm{NS}^G(M)_\mathbb{R}$ is $G$-big if and only if
$\xi=a+b$ for some
    $a\in \mathrm{Int}\left
(\mathrm{C}^G(M)\right )$ and $b\in \mathrm{Eff}^G(M)$.
\end{prop}

This leads to the following analogue of Theorem 2.2.26 of
\cite{lazI} (the proof is similar):


\begin{thm}\label{thm:big-and-eff-closure}
Let $M$ be a complex projective $G$-manifold, and suppose that there
exists an ample $B\in \mathrm{Pic}^G(M)$ such that
$M^\mathrm{s}(B)\neq \emptyset$, $\mathrm{codim}\left
(M^\mathrm{u}(B)\right )\ge 2$. Then the $G$-big cone is open.
Furthermore, $\mathrm{Big}^G(M)$ is the interior of the
$G$-pseudoeffective cone, and the $G$-pseudoeffective cone is the
closure of the $G$-big cone:
$$
\mathrm{Big}^G(M)=\mathrm{Int}\left
(\overline{\mathrm{Eff}^G(M)}\right ),\,\,\,\,\,\,\,\,\,
\overline{\mathrm{Eff}^G(M)}=\overline{\mathrm{Big}^G(M)}.$$
\end{thm}


\begin{rem}\label{rem:generically-free-actions-vol-mu-extends}
Let us take up again the special case of generically free actions,
under hypothesis of Corollary \ref{cor:vol-mu-pos-equiv-charat}. By
Corollary \ref{cor:big-and-locally-free-homogeneous} and the
argument used in Remark \ref{rem:well-defined-over-Q}, for every
$\mu$ the function $\mathrm{vol}_\mu$ extends to a well defined
homogeneous  function of degree $\mathrm{n-g}$ on
 the set $\mathrm{Big}^G(M)_{\mathbb{Q}}\subseteq \mathrm{Big}^G(M)$
 of rational points in the
$G$-big cone. We shall denote this function by
$\mathrm{vol}_\mu:\mathrm{Big}^G(M)_{\mathbb{Q}}\rightarrow
\mathbb{R}$.
\end{rem}


\section{Continuity}

We have noted in Remark \ref{rem:well-defined-over-Q} that the
functions $\mathrm{vol}_\mu:\mathrm{Pic}^G(M)\rightarrow \mathbb{R}$
descend to well-defined functions on $\mathrm{Num}^G(M)$, and that
by homogeneity $\mathrm{vol}_0$ extends to a function $
\mathrm{vol}_0:\mathrm{Num}^G(M)_\mathbb{Q}\rightarrow \mathbb{R}$
homogeneous of degree $\mathrm{n-g}$. We have also noted in Remark
\ref{rem:generically-free-actions-vol-mu-extends} that, if the
action of $G$ on $M$ is generically free, then $\mathrm{vol}_\mu$
extends to a function on $\mathrm{Big}^G(M)_\mathbb{Q}$ for every
$\mu$, also homogeneous of degree $\mathrm{n-g}$.

We shall now dwell on the continuity properties of these functions
and prove the estimates (\ref{eqn:continuity-bound-trivial-rep}) and
(\ref{eqn:continuity-bound-general-rep}) below. We shall deduce
that, at least under appropriate hypothesis, $\mathrm{vol}_0$
extends to a continuous real-valued function on $
\mathrm{NS}^G(M)_\mathbb{R}=:\mathrm{Num}^G(M)_\mathbb{Q}\otimes
\mathbb{R}$, and that $\mathrm{vol}_\mu$ extends to a continuous
real-valued function on $ \mathrm{Big}^G(M)\subseteq
\mathrm{NS}^G(M)_\mathbb{R}$.

With some abuse, we shall not distinguish notationally an element
$A\in \mathrm{Pic}^G(M)$ from its equivalence class in
$\mathrm{Num}^G(M)_\mathbb{Q}\subseteq \mathrm{NS}^G(M)_\mathbb{R}$.

The content of this section is given by the following two Theorems:


\begin{thm}\label{thm:continuity}
Let $M$ be a complex projective $G$-manifold. Suppose that there
exists $B\in \mathrm{Pic}^G(M)$ ample and regular, satisfying
$\mathrm{codim}\big (M^\mathrm{u}(B)\big )\ge 2$. Let us fix any
norm on $\mathrm{NS}^G(M)_\mathbb{R}$. Then there exists $C>0$ such
that for any $D,D'\in \mathrm{NS}^G(M)_{\mathbb{Q}}$ we have:
\begin{eqnarray}
\label{eqn:continuity-bound-trivial-rep} \left |
\mathrm{vol}_0(D)-\mathrm{vol}_0(D')\right |\le C\, \max\{
\|D\|,\|D'\|\}^{\mathrm{n-g}-1}\cdot \|D-D'\|.
\end{eqnarray}

\end{thm}


\begin{thm}\label{thm:continuity-general-rep}
In the situation of Theorem \ref{thm:continuity}, assume in addition
that the action of $G$ on $M$ is generically free. Then there exists
$C>0$ such that for any $D,D'\in \mathrm{Big}^G(M)_{\mathbb{Q}}$ and
every highest weight $\mu$ we have:
\begin{eqnarray}
\label{eqn:continuity-bound-general-rep} \left |
\mathrm{vol}_\mu(D)-\mathrm{vol}_\mu(D')\right |\le C\, \dim
(V_\mu)^2\,\max\{ \|D\|,\|D'\|\}^{\mathrm{n-g}-1}\cdot \|D-D'\|.
\end{eqnarray}

\end{thm}


We shall give an explicit proof of Theorem \ref{thm:continuity}, and
leave it to the reader to make the necessary changes for the proof
of Theorem \ref{thm:continuity-general-rep}. The statement and the
proof of Theorem \ref{thm:continuity} are inspired by Theorem 2.2.44
of \cite{lazI}. We shall use throughout additive notation.

\textit{Proof.} By assumption, there exists an open (conic) chamber
$C$ in the $G$-ample cone $C^G(M)\subseteq
\mathrm{NS}^G(M)_\mathbb{R}$ such that $B\in C$ and every $A\in
\mathrm{Pic}^G(M)$ whose numerical class is in $C$ satisfies
$M^\mathrm{s}(A)=M^\mathrm{ss}(A)=M^\mathrm{s}(B)\neq \emptyset$.
Setting $r=:\dim \mathrm{NS}^G(M)_\mathbb{R}$, we can then find
$A_1,\ldots,A_r\in \mathrm{Pic}^G(M)$ whose numerical classes all
belong to $C$, and such that $\mathcal{B}=\{A_1,\ldots, A_r\}$ is a
basis of $\mathrm{NS}^G(M)_\mathbb{R}$. Perhaps after replacing each
$A_j$ by a sufficiently large multiple we may, and shall, suppose
that the $A_j$'s are all very ample and descend to very ample line
bundles on the GIT quotient $M_0(A_j)=M_0(B)$.

As in \cite{lazI}, we shall use the norm $\|\xi\|=\max \{|x_i|\}$,
where $\xi =\sum_ix_i\cdot A_i$.

\begin{prop}\label{prop:integral-linear-combis}
There exists $C>0$ such that for any $D,D'\in \mathrm{NS}^G(M)$ of
the form
\begin{eqnarray*}
D&=&a_1\cdot A_1+\ldots+a_r\cdot A_r\\
D'&=&a'_1\cdot A_1+\ldots+a'_r\cdot A_r,
\end{eqnarray*}
where $a_i,a_i'\in \mathbb{Z}$ for every $i$, we have:
\begin{eqnarray}
\label{eqn:continuity-bound-general-rep-integer-combi} \left |
\mathrm{vol}_0(D)-\mathrm{vol}_0(D')\right |\le C\, \max\{
\|D\|,\|D'\|\}^{\mathrm{n-g}-1}\cdot \|D-D'\|.
\end{eqnarray}
\end{prop}

\textit{Proof of Proposition \ref{prop:integral-linear-combis}.} By
abuse of language, it will be convenient to identify each $A_j$ with
some general divisor in the linear series $|A_j|^G$. With this
interpretation, every $A_j$ is a $G$-invariantly effective divisor
in $M$, non-singular away from $M^\mathrm{u}(B)$. We shall assume to
begin with that $a_i\ge a'_i$ for every $i$; thus $b_i=:a_i-a'_i\ge
0$ for every $i$ and $D'=D-B$, where $B=\sum _ib_i\cdot A_i$ is
$G$-invariantly effective. Thus, $\mathrm{vol}_0(D-B)\le
\mathrm{vol}_0(D)$.

In order to obtain a bound in the opposite direction, let us choose
for every $j$ a very general $E_j\in |A_j|^G$. Then $E_j$ is
$G$-invariantly effective, reduced, and meets every $A_i$ properly.
Let us set $D_j=\left .D\right |_{E_j}$. By Proposition
\ref{prop:invariant-hypersurfaces-finite-vol-codimension-two},
$\mathrm{vol}_0(E_j,D_j)<+\infty$.

\begin{lem}\label{lem:first-lemma-on-continuity-bound}
Given that $b_i\ge 0$ for every $i$, we have
\begin{equation}\label{eqn:first-bound-D-B}
\mathrm{vol}_0(D-B)\ge \mathrm{vol}_0(D)-(n-g)\cdot \sum _jb_j \cdot
\mathrm{vol}_0 (E_j,D_j).\end{equation}
\end{lem}

\textit{Proof of Lemma \ref{lem:first-lemma-on-continuity-bound}.}
It suffices by induction to consider the case $B=b_1A_1$.

For $m\gg 0$, let us choose $mb_1$ very general divisors $F_\alpha
\in |A_1|^G$. This yields an equivariant exact sequence
$$
0\rightarrow \mathcal{O}_M(mD-mb_1A_1)\rightarrow
\mathcal{O}_M(mD)\rightarrow \bigoplus _{\alpha
=1}^{mb_1}\mathcal{O}_{F_\alpha}(mD).$$ By the very generality of
our choices, we have an estimate
\begin{eqnarray*}
\dim H^0\big (F_\alpha,\mathcal{O}_{F_\alpha}(mD)\big )^G&=&\dim
H^0\big (E_1,\mathcal{O}_{E_1}(m D_1)\big )^G\\
&\le &\mathrm{vol}_0 (E_1,D_1)\cdot
\frac{m^{\mathrm{n-g}-1}}{(\mathrm{n-g}-1)!}+o(m^{\mathrm{n-g}-1}).\end{eqnarray*}
Thus, setting $h^0(M,L)_0=:\dim H^0(M,L)^G$ for $L\in
\mathrm{Pic}^G(M)$,
\begin{eqnarray*}
\frac{(\mathrm{n-g})!}{m^{\mathrm{n-g}}}\,h^0 \big
(M,m(D-b_1A_1)\big
)_0&\ge&\frac{(\mathrm{n-g})!}{m^{\mathrm{n-g}}}\,h^0 \big (M,mD\big
)_0\\
&&-b_1\,(\mathrm{n-g})\,\frac{(\mathrm{n-g}-1)!}{m^{\mathrm{n-g-1}}}\,h^0
\big (E_1,D_1\big )_0\\
&\ge&\frac{(\mathrm{n-g})!}{m^{\mathrm{n-g}}}\,h^0 \big (M,mD\big
)_0\\
&&-(\mathrm{n-g})\,b_1\,\mathrm{vol}_0(E_1,D_1)+o(1).
\end{eqnarray*}
and taking $\limsup$ as $m\rightarrow +\infty$ we obtain:
\begin{eqnarray*}
\mathrm{vol}_0 \big (\mathcal{O}_M(D-b_1A_1)\big
)&\ge&\mathrm{vol}_0 \big (\mathcal{O}_M(D)\big
)-b_1\,(\mathrm{n-g})\,\mathrm{vol}_0 \big (E_1,D_1\big ).
\end{eqnarray*}
This completes the proof of Lemma
\ref{lem:first-lemma-on-continuity-bound}.

\begin{prop}\label{prop:equivariant-bound-norm}
There exists $C>0$ such that
$$
\mathrm{vol}_0\big (E_j,D_j\big )\le C\, \|D\|^{\mathrm{n-g}-1}$$
for every $j=1,\ldots,r$.
\end{prop}

\textit{Proof of Proposition \ref{prop:equivariant-bound-norm}.} Let
$E_j=\sum_{k=1}^{r_j}E_{jk}$ be the decomposition of $E_j$ into its
irreducible components, and for every $j$ and $k$ let us set
$D_{jk}=:\left .D\right |_{E_{jk}}$ (notice that $E_j$ may be
assumed irreducible as soon as $\mathrm{n-g}>1$). Since
$\mathrm{vol}_\mu\big (E_j,D_j\big )\le
\sum_{k=1}^{r_j}\mathrm{vol}_\mu\big (E_{jk},D_{jk}\big )$, it
suffices to prove the statement for each irreducible component
$E_{jk}$. Let us set
$$D^+=:\sum _{i=1}^r |a_i|D_i,\,\,\,\,\,D^+_j=:\left .D^+\right
|_{E_j}.$$ Then
$$D^+_{jk}-D_{jk}=\sum _i (|a_i|-a_i)\left.A_i\right |_{E_{jk}},$$ where each $
\left.A_i\right |_{E_{jk}}$ is a $G$-invariantly effective line
bundle on $E_{jk}$. It follows by Remark \ref{rem:G-eff-morphism}
that
$$
\mathrm{vol}_0\big ( E_{jk},D_{jk} \big ) \le \mathrm{vol}_0\big (
E_{jk},D_{jk}^+\big ).
$$
Since $\|D\|=\|D^+\|$, it thus suffices to prove the statement of
Proposition \ref{prop:equivariant-bound-norm} using $D^+$ in place
of $D$. To this end, let us first record the following:

\begin{lem} \label{lem:auxiliary-transl}
There exists $R\in \mathrm{GL}\big ( \mathrm{NS}^G(M)_\mathbb{R}\big
)$ with the following properties:
\begin{description}
  \item[i):] for every $i=1,\ldots,r$, $R(A_i)-A_i$ represents
  a $G$-invariantly effective positive power of $B$;
  \item[ii):] for every $i=1,\ldots,r$, $R(A_i)$ lies in the
  interior of the chamber $C$.
\end{description}
In particular, $R(D^+)$ also lies in the
  interior of the chamber $C$.
\end{lem}

\textit{Proof of Lemma \ref{lem:auxiliary-transl}.} We need only
define $R(A_i)=:A_i+re_G(B)B$ for some $r\gg 0$, and then extend by
linearity. This establishes i), and ii) is a consequence of i) given
our choice of the $A_i$'s.

\bigskip

Back to the proof of Proposition \ref{prop:equivariant-bound-norm},
it follows from Lemma \ref{lem:auxiliary-transl} that $R(D^+)-D^+$
is $G$-invariantly effective, and that so is its restriction to
$E_j$ for every $j$. Therefore,
\begin{equation}\label{eqn:further-reduction-with-R}
\mathrm{vol}_0\big ( E_{jk},D_{jk} \big ) \le \mathrm{vol}_0\big (
E_{jk},D_{jk}^+\big )\le \mathrm{vol}_0\big ( E_{jk},\left
.R(D^+)\right |_{E_{jk}} \big ).
\end{equation}

To fix ideas, let us first suppose that all the $E_{jk}$'s are
non-singular. If $M_0(B)$ is the GIT quotient of $M$ with respect to
the linearization $A_j$, the GIT quotient of $E_{jk}\subseteq M$
with respect to the same linearization is an ample divisor
$(E_{jk})_0\subseteq M_0(B)$; hence it is $(
\mathrm{n}-\mathrm{g}-1)$-dimensional .

Since $\left .R(D^+)\right |_{E_{jk}}$ is an ample and regular
$G$-linearized line bundle, by Corollary
\ref{cor:main-regular-case-but-general} and Remark
\ref{rem:main-regular-case-but-general} we have
\begin{eqnarray*}
\mathrm{vol}_0\left ( E_{jk},\left .R(D^+)\right |_{E_{jk}} \right )
& =& \mathrm{vol}_0\big ( (E_{jk})_0,\Omega
_{R(D^+)}^{(jk)}\big )\nonumber \\
&=&C'_{jk}\,\int _{(E_{jk})_0}\big (\Omega _{R(D^+)}^{(jk)}\big )
^{\mathrm{n-g}-1} \\
&\le& C''_{jk}\,\|D^+\|^{\mathrm{n-g}-1},\end{eqnarray*}
where $ \Omega _{R(D^+)}^{(jk)}$ is a K\"{a}hler form on $(E_j)_0$,
which we may choose to depend linearly on $R(D^+)$, whence on $D^+$
($ \Omega _{R(D^+)}^{(jk)}$ is the restriction to
$(E_{jk})_0\subseteq M_0(B)$ of a K\"{a}hler form on $M_0(B)$ of the
form $\sum _i\ell_i\Omega _{A_i}$, where $R(D^+)=\sum _i\ell
_iA_i$).

Let us now argue in general.

\begin{lem}\label{lem:generality-of-Ej}
Given the generality in its choice, $E_{jk}$ satisfies the following
properties:
\begin{description}
  \item[i):] $E_{jk}$ is non-singular away from
$M^{\mathrm{u}}(A_j)=M^{\mathrm{u}}(R(D^+))=M^\mathrm{u}(B)$;
  \item[ii):] if $\Phi_{R(D^+)}:M\rightarrow \mathfrak{g}^*$ is the
moment map of $R(D^+)$, $E_{jk}$ is non-singular in a neighbourhood
of $E_{jk}\cap \Phi_{R(D^+)}^{-1}(0)$;
  \item[iii):] $E_{jk}$ is transversal
to $\Phi_{R(D^+)}^{-1}(0)$.
\end{description}
\end{lem}

\textit{Proof of Lemma \ref{lem:generality-of-Ej}.} i) follows
directly from Bertini's Theorem and the definition of unstable
locus.

Given i), ii) follows from the well-known relation between the
semi-stable locus and the zero locus of the moment map \cite{kir1}.

As to iii), by the arguments of Lemma 3 in \cite{pao-mm} and
compactness one can see the following: There exist a finite number
of holomorphic embeddings  $\varphi
_i:B_\mathrm{n-g}(0,1)\rightarrow M$, where
$B_\mathrm{n-g}(0,1)\subseteq \mathbb{C}^{\mathrm{n-g}}$ is the unit
ball centered at the origin, satisfying a): $\varphi _i\left
(B_\mathrm{n-g}(0,1)\right)\subseteq \Phi_{R(D^+)}^{-1}(0)$; b): as
a submanifold of $\Phi_{R(D^+)}^{-1}(0)$, $\varphi
_i(B_\mathrm{n-g}(0,1))$ is transversal to every $G$-orbit; c): the
union $\bigcup _i \varphi _i(B_\mathrm{n-g}(0,1))$ maps surjectively
onto $M_0\big (R(D^+)\big)= M_0(B)$. In view the local analytic
proof of Bertini's theorem in \cite{gh}, we may assume that $E_j$ is
transversal to each $\varphi _i(B)$. By $G$-invariance, it is then
transversal to all of $\Phi_{R(D^+)}^{-1}(0)$.

This completes the proof of Lemma \ref{lem:generality-of-Ej}.

\bigskip

Let $f_{jk}:\widetilde{E}_{jk}\rightarrow E_{jk}$ be a
$G$-equivariant resolution of singularities, which may be assumed to
be an isomorphism away from $M^\mathrm{u}\big (R(D^+)\big )=
M^\mathrm{u}(B)$.

\begin{cor}\label{cor:generality-of-Ej}
$0\in \frak{g}^*$ is a regular value of the moment map for
$f_{jk}^*\big (R(D^+)\big )$ given by
$$\widetilde{\Phi}_{R(D^+)}^{(jk)}=:\Phi_{R(D^+)}\circ f_{jk} :
\widetilde{E}_{jk}\longrightarrow \frak{g}^*,$$ and $\left (
\widetilde{\Phi}_{R(D^+)}^{(jk)}\right)^{-1}(0)\neq \emptyset$.
Furthermore, $\widetilde{\Phi}_{R(D^+)}^{(jk)}$ is bounded in norm
away from zero on $f_{jk}^{-1}\big (M^\mathrm{u}(B)\big )$.
\end{cor}

\textit{Proof of Corollary \ref{cor:generality-of-Ej}.} The only
statement not immediately obvious from Lemma
\ref{lem:generality-of-Ej} is that $\left (
\widetilde{\Phi}_{R(D^+)}^{(jk)}\right)^{-1}(0)\neq \emptyset$. To
see this, we need only recall that $$M^{s}\left(\left. R(D^+)\right
|_{E_{jk}}\right )=M^{\mathrm{s}}\left(\left.B\right
|_{E_{jk}}\right )\neq \emptyset.$$

Returning to the proof of Proposition
\ref{prop:equivariant-bound-norm}, there exists a $G$-invariantly
effective exceptional divisor $F_{jk}\subseteq \widetilde{E}_{jk}$
such that $f_j^*(R(D^+))^{\otimes l}(-F_{jk})\in \mathrm{Pic}^G(
\widetilde{E}_{jk})$ is a (very) ample $G$-linearized line bundle on
$\widetilde{E}_{jk}$, for every $l\gg 0$. Its linearization
corresponds to a moment map (with respect to an appropriate
K\"{a}hler form) of the form
$$\widetilde{\Phi} _{ljk}=:l\,\widetilde{\Phi}_{R(D^+)}^{(j)}+\Phi _{F_{jk}}:
\widetilde{E}_{jk}\rightarrow \frak{g}^*,$$ where $\Phi _{F_{jk}}$
is an appropriate fixed equivariant map. By the above, if $l\gg 0$
then $0\in \frak{g}^*$ is a regular value of $\widetilde{\Phi}
_{ljk}$ and $(\widetilde{\Phi} _{ljk})^{-1}(0)\neq \emptyset$. Thus,

\begin{lem}\label{lem:ample-and-reg-tilde}
If $l\gg 0$, $f_j^*\big (R_j(D^+)\big )^{\otimes l}(-F_{jk})\in
\mathrm{Pic}^G(\widetilde{E}_{jk})$ is ample and regular.\end{lem}

Let us then fix $l_0\gg 0$ such that $f_{jk}^*\big (R(D^+)\big
)^{\otimes l_0}(-F_{jk})$ is ample and regular, and next $a_0\gg 0$
such that in addition $f_{jk}^*\big (R(D^+)\big )^{\otimes
a_0l_0}(-a_0F_{jk})$ is $G$-invariantly effective.

Next let $p_i\uparrow +\infty$ be a sequence of integers. In view of
Corollary \ref{cor:general-homog-trivial-rep} we have
$\mathrm{vol}_0 \left (\left . R(D^+)\right |_{E_{jk}}^{\otimes
p_i}\right )=p_i^{\mathrm{n-g}-1}\mathrm{vol}_0\big (\left .
R(D^+)\right |_{E_j}\big) $ for every $i$. On the upshot, we have
\begin{eqnarray}\label{eqn:equivariant-bound-norm-that's-it}
\mathrm{vol}_0\big (\left . R(D^+)\right |_{E_{jk}}\big)
&\le&\mathrm{vol}_0\left (f^*_{jk} \big(R(D^+)\big)\right) \\
&=&\lim _{i\rightarrow +\infty}\frac{\mathrm{vol}_0 \left
(f_{jk}^*\big(R(D^+)\big )^{\otimes p_i}\right
)}{p_i^{\mathrm{n-g}-1}}\nonumber\\
&\le&\limsup _{i\rightarrow +\infty}\frac{\mathrm{vol}_0 \left
(f_{jk}^*\big(R(D^+)\big )^{\otimes (p_i+a_0l_0)}(-a_0F_{jk})\right
)}{p_i^{\mathrm{n-g}-1}}.\nonumber
\end{eqnarray}
In view of Lemma \ref{lem:ample-and-reg-tilde}, the volume in the
last line of (\ref{eqn:equivariant-bound-norm-that's-it}) may be
computed using Theorem \ref{thm:main-regular-case}. The K\"{a}hler
form for $ f_{jk}^*\big(R(D^+)\big )^{\otimes
(p_i+a_0l_0)}(-a_0F_{jk})$ may clearly be chosen of the form $
(p_i+a_0l_0)\Omega_{R(D^+)}+\Omega _F$, where $\Omega_{R(D^+)}$
depends linearly on $R(D^+)$, and $\Omega_F$ is fixed. Thus, given
(\ref{eqn:equivariant-bound-norm-that's-it}), we obtain
\begin{eqnarray*}
\mathrm{vol}_0\big (\left . R(D^+)\right
|_{E_{jk}}\big)&\le&\lim_{i\rightarrow +\infty} \int _{\widetilde{
E}_{jk}}\frac{\left ((p_i+a_0l_0)\Omega_{R(D^+)}+\Omega
_F\right ) ^{\mathrm{n-g}-1}}{p_i^{\mathrm{n-g}-1}}\\
&\le &C\|D^+\|^{\mathrm{n-g}-1}.\end{eqnarray*} Since the $E_{jk}$'s
are in finitely many, the constant $C$ may be chosen independent of
$j$ and $k$. The statement of Proposition
\ref{prop:equivariant-bound-norm} follows.

\bigskip
On the upshot, there exists a constant $C>0$ such that whenever
$D=\sum _ia_iA_i$ with $a_i\in \mathbb{Z}$ for every $i$, and
$B=\sum b_iA_i$ with $b_i\in \mathbb{N}$ for every $i$, we have
\begin{equation}
\label{eqn:first-bound-effective-case}
\mathrm{vol}_0(D)-C\,\|D\|^{\mathrm{n-g}-1}\cdot \|B\|\le
\mathrm{vol}_0(D-B)\le \mathrm{vol}_0(D). \end{equation}

Given this, we can now prove
(\ref{eqn:continuity-bound-general-rep-integer-combi}) without
assuming that $B=D-D'$ is effective, as follows. Let us write
$B=F-E$, where $F=\sum _if_iA_i$, $E=\sum_ie_iA_i$ and $f_i,e_i\in
\mathbb{N}$ for every $i$. We may also assume that $e_if_i=0$ for
every $i$. Then (\ref{eqn:first-bound-effective-case}) implies the
inequalities:
\begin{eqnarray*}
\mathrm{vol}_0(D)-\mathrm{vol}_0(D-F)&\le& C\,\|D\|^{\mathrm{n-g}-1}\cdot \|F\|\\
\mathrm{vol}_0(D+E-F)-\mathrm{vol}_0(D-F)&\le&
C\,\|D+E-F\|^{\mathrm{n-g}-1}\cdot \|E\|.\end{eqnarray*} The
statement then follows by the triangle inequality since $\max \left
(\|E\|,\|F\|\right)=\|E-F\|$. This completes the proof of
Proposition \ref{prop:integral-linear-combis}.

\bigskip

Back to the proof of Theorem \ref{thm:continuity}, now we need only
notice that the terms in
(\ref{eqn:continuity-bound-general-rep-integer-combi}) are all
homogeneous of degree $\mathrm{n-g}$. Therefore, the inequality in
the statement of Theorem \ref{thm:continuity} holds for any $D,D'\in
\mathrm{NS}^G(M)_\mathbb{Q}$.


\begin{cor} \label{cor:extension-by-continuity}
In the hypothesis of Theorem \ref{thm:continuity}, $\mathrm{vol}_0$
extends to a continuous real-valued function on $
\mathrm{NS}^G(M)_\mathbb{R}\setminus \{0\}$, and on $
\mathrm{NS}^G(M)_\mathbb{R}$ if $\mathrm{n\ge g}+1$.\end{cor}


\begin{cor} \label{cor:extension-by-continuity-general-repr}
In the hypothesis of Theorem \ref{thm:continuity-general-rep},
$\mathrm{vol}_\mu$ extends to a continuous real-valued function on $
\overline{\mathrm{Big}^G(M)}\setminus \{0\}$, and on $
\overline{\mathrm{Big}^G(M)}$ if $\mathrm{n\ge g}+1$.\end{cor}

\section{Equivariant Fujita approximations}

In the study of absolute algebro-geometric volumes, an important
role is played by Fujita approximations, which can be seen as
providing a sort of asymptotic Zariski decomposition of an arbitrary
big class \cite{f}, \cite{del}, \cite{lazI}. In this section, we
shall give an equivariant version of this in the special case of
circle actions, under some assumptions on the action of $S^1$ on
$M$. An equivariant version of Fujita approximations has been given
in \cite{pao-ages} for finite group actions.

\begin{thm}\label{thm:equiv-fujita-approx}
Let $M$ be a projective $S^1$-manifold, and suppose that there
exists a regular and ample $B\in \mathrm{Pic}^{S^1}(M)$ such that
$\mathrm{codim}\big (M^\mathrm{u}(B)\big )\ge 2$. Suppose $L\in
\mathrm{Pic}^{S^1}(M)$ is $S^1$-big, and fix $\epsilon >0$. Let
$F\subseteq \mathbb{Z}$ be any finite set of weights, with $0\in F$.
Then there exist:
\begin{itemize}
  \item an
equivariant modification $\beta :\widetilde{M}\rightarrow M$
(depending on $\xi$, $F$ and $\epsilon$), where the $G$-ample cone
of $\widetilde{M}$ has non-empty interior,
  \item $A,E\in \mathrm{Pic}^{S^1}(\widetilde{M})$ with $A$ ample
  and $E$ $S^1$-invariantly
  effective,
  \item an integer $p>0$ with $(p,e_G(L))=1$,
\end{itemize}
such that $ \beta ^*(L)^{\otimes p}=A\otimes E$,
$\mathrm{vol}_\mu(A)\ge p^{\mathrm{n-1}}\big (
\mathrm{vol}_\mu(L)-\epsilon\big )$ for every $\mu \in F$, and
$\mathrm{vol}(A)\ge p^{\mathrm{n}}\big (
\mathrm{vol}(L)-\epsilon\big )$.
\end{thm}

\textit{Proof.} We shall subdivide the proof in a series of
Propositions.

\begin{prop}\label{prop:fujita-step-1}
A decomposition as in the statement of Theorem
\ref{thm:equiv-fujita-approx} exists, if we only require $A\in
\mathrm{Pic}^{S^1}(\widetilde{M})$ to be nef, big and $S^1$-big.
\end{prop}

\textit{Proof.} The argument below is inspired by the proof of
Theorem 11.4.4 of \cite{lazI}.

Perhaps after replacing $B$ by some sufficiently large and divisible
power, we may assume that $e_{S^1}(B)=1$, $B$ is very ample and $
K_M\otimes B^{\otimes (\mathrm{n}+1)}\in \mathrm{Pic}^{S^1}(M)$ is
ample and GIT-equivalent to $B$. Let us then define $S=:\Big
(K_M\otimes B^{\otimes (\mathrm{n}+1)}\Big ) ^{\otimes e}$, where
$e=:e_{S^1}\big (K_M\big )$. By arguing as in the proof of Lemma
\ref{lem:kod-auxiliary}, we may also suppose without loss of
generality that $H^0\big (M,S\big)^{S^1}\neq \{0\}$. Let us choose
$\sigma \in H^0\big (M,S\big)^{S^1}\neq \{0\}$, $\sigma \neq 0$.

For $p\ge 0$, let $R_p=:L^{\otimes p}\otimes S^{-1}\in
\mathrm{Pic}^{S^1}(M)$. By Proposition
\ref{prop:limit-infinity-p-prime} and Remark
\ref{rem:lim-sup-and-lim}, there exists $p\gg 0$ prime with
$e_{S^1}(L)$ such that $\mathrm{vol}_\mu(R_p)>p^{\mathrm{n-g}}\big (
\mathrm{vol}_\mu(L)-\epsilon\big )$ for every $\mu \in F$. In
particular, $R_p$ is $S^1$-big; by Corollary
\ref{cor_G-big-implies-big}, it is big.

Let now $ \mathcal{J}=\mathcal{J}(\|R_p\|)$ be the asymptotic
multiplier ideal of $R_p$ \cite{lazI}. Since $R_p\in
\mathrm{Pic}^{S^1}(M)$, $ \mathcal{J}\subseteq \mathcal{O}_M$ is
$G$-invariant. Let $\beta:\widetilde{M}\rightarrow M$ be a
resolution of the blow-up of $ \mathcal{J}$; we may assume that
$\mu$ is given by a finite sequence of blow-ups along
$S^1$-invariant smooth centers, and therefore it is an equivariant
birational morphism of projective $S^1$-manifolds. Thus
$\mathcal{J}\cdot
\mathcal{O}_{\widetilde{M}}=\mathcal{O}_{\widetilde{M}}(-E_p)$, for
an appropriate $S^1$-invariantly effective divisor $E_p\subseteq
\widetilde{M}$. Set $A_p=:\beta ^*(L) ^{\otimes p}(-E_p)\in
\mathrm{Pic}^{S^1}\big ( \widetilde{M}\big )$.

\begin{lem}\label{lem:A-is-nef-on-log-res}
$A_p$ is globally generated (hence nef), big and $S^1$-big.
\end{lem}

\textit{Proof.} Since, in additive notation, $pL=K_M+R_p+nB+\Big \{
(e-1)\big [K_X+(n+1)B\big ]+B\Big\}$, Corollary 11.2.13 of
\cite{lazI} implies that $L^{\otimes p}\otimes \mathcal{J}$ is
globally generated. This in turn implies that $A_p$ is globally
generated.

\begin{lem}\label{lem:equiv-emb-fujita}
For every $\ell \ge 1$, there is an equivariant injective linear map
$H^0\big ( M,R_p^{\otimes \ell}\big )\rightarrow H^0\big (
M,L^{\otimes p\ell}\otimes \mathcal{J}^\ell\big )$.\end{lem}

\textit{Proof.} This follows from the same arguments as in
\cite{lazI}, because all the sheaves involved and the section
$\sigma$ are $S^1$-invariant.

\bigskip

As in \cite{lazI}, one can deduce from this an equivariant injective
linear map $H^0\big (M,R_p^{\otimes \ell}\big )\rightarrow H^0\big
(\widetilde{M},A_p^{\otimes \ell}\big )$. This implies the statement
of Proposition \ref{prop:fujita-step-1} with $A=A_p$,
$E=\mathcal{O}_{\widetilde{M}}(E_p)$.

\hfill Q.E.D.

\bigskip

Since when passing to a resolution one might lose control of the
codimension of unstable loci, we may not resort directly to the
equivariant Kodaira Lemma to complete the proof of the Theorem by
analogy with the action-free case. We shall therefore adopt an
\textit{ad hoc} argument.

\begin{prop}\label{prop:fujita-step-2}
A decomposition as in the statement of Theorem
\ref{thm:equiv-fujita-approx} exists, if we only require $A\in
\mathrm{Pic}^{S^1}(\widetilde{M})$ to be ample, big and $S^1$-big.
\end{prop}

\textit{Proof.} Let $\widetilde{M}$, $A=A_p$, $E$, $p$ be as in
Proposition \ref{prop:fujita-step-1}. Let $\varphi
_k:\widetilde{M}\rightarrow \mathbb{P}H^0(\widetilde{M},A_p^{\otimes
k})^*$ be the $S^1$-equivariant projective morphism induced by the
linear series $\big |A_p^{\otimes k}\big |$, $k\ge 1$. Let
$\widetilde{M}_k=:\varphi _k\big (M\big )$,
$H_k=\mathcal{O}_{M_k}(1)$. For $k\gg 0$, the projective morphisms
$\varphi _k:\widetilde{M}\rightarrow \widetilde{M}_k$ stabilize to
an $S^1$-equivariant algebraic fibre space $\varphi
_\infty:\widetilde{M}\rightarrow \widetilde{M}_\infty$, and $\varphi
_k^*\big ( H_k\big )=A^{\otimes k}$ for every $k\ge 0$. One can
deduce that there exists $H\in \mathrm{Pic}^{S^1}\big
(\widetilde{M}_\infty\big )$ ample and such that $\varphi _\infty
^{*}\big (H\big )=A$.

By composing a finite sequence of blow-ups along $S^1$-invariant
smooth centers, we can find an $S^1$-equivariant resolution of
singularities of the inverse birational map $\varphi _\infty^{-1}$,
$\widetilde{M}\stackrel{\psi}{\leftarrow} \widehat{M}_\infty
\stackrel{\phi}{\rightarrow}\widetilde{M}_\infty$. In particular,
$\phi =\phi _\infty \circ \psi$. It follows that there exists an
$S^1$-invariantly effective divisor $F\subseteq
\widetilde{M}_\infty$ such that $\psi ^*\big (A^{\otimes k}\big
)(-F)\in \mathrm{Pic}^{S^1}( \widetilde{M}_\infty)$ is ample for
every $k\gg 0$. Suppose that $F=\sum _ia_iF_i$, where each $a_i>0$
and $F_i\subseteq F$ is irreducible.

\begin{lem}\label{lem:exists-a-section}
We can find $k_0$ (divisible by $e_{S^1}(A)=e_{S^1}(L)$) such that
$$H^0\left (\widehat{M}_\infty,\psi ^*\big ( A_p\big )^{\otimes
k_0}(-F)\right )^{S^1}\neq 0.$$
\end{lem}

\textit{Proof of Lemma \ref{lem:exists-a-section}.} By taking tensor
products, we are reduced to proving the statement in the case where
$F$ is reduced and irreducible, i.e., one of the $F_i$'s. Since $F$
is $S^1$-invariantly effective, we have for every $l\ge 0$ an
equivariant short exact sequence of sheaves on $\widehat{M}_\infty$,
\begin{equation}\label{eqn:exact-on-hat-M}
0\rightarrow \psi ^*\big ( A\big )^{\otimes l}(-F)\longrightarrow
\psi ^*\big ( A\big )^{\otimes l}\longrightarrow \psi ^*\big ( A\big
)^{\otimes l}\otimes \mathcal{O}_F\rightarrow 0.\end{equation} Now $
\mathrm{vol}_0\big(\psi ^* ( A )\big )= \mathrm{vol}_0(A)>0$, while:

\begin{claim}\label{claim:kay-estimate-n-g-1}
For some constant $C>0$, we have $$\dim H^0\Big (F,\psi ^*(
A)^{\otimes k}\otimes \mathcal{O}_F\Big )^G\le C\,k^{\mathrm{n}-2}$$
for $k\gg 0$.\end{claim}

\textit{Proof.} We have $\psi ^*(A)=\phi ^*(H)$. Let $S=:\psi
(F)\subseteq \widetilde{M}_\infty$, $\gamma:F\rightarrow S$ the
induced morphism. Thus $S$ is an $S^1$-invariant irreducible
subvariety of $\widetilde{M}_\infty$, and $\dim (S)\le
\mathrm{n}-2$. By the projection formula and a free resolution of
$\gamma_*\mathcal{O}_F$ by powers of $H$, we are reduced to proving
that $\dim H^0\big(S,H^{k}\otimes \mathcal{O}_S\big )^{S^1}\le
C\,k^{\mathrm{n-2}}$, and this is obvious by dimension reasons.

\bigskip

The statement of Lemma \ref{lem:exists-a-section} then follows by
taking global sections in (\ref{eqn:exact-on-hat-M}).

\hfill Q.E.D.

\bigskip

Now let us set $ \widetilde{\beta}=:\beta \circ \psi:
\widehat{M}_\infty\rightarrow M$. Then, for every $k\ge 0$ we have:
\begin{eqnarray*}
\widetilde{\beta}^*\big (L\big ) ^{\otimes p^k}&=&\psi
^*(A)^{\otimes p^{k-1}}\otimes \psi ^*(E)^{\otimes p^{k-1}}\\
&=&\psi ^*(A)^{\otimes (p^{k-1}-k_0)}\otimes \psi ^*(A)^{\otimes
k_0}(-F)\otimes
\Big (F\otimes \psi ^*(E)^{\otimes p^{^{k-1}}}\Big )\\
&=&\widetilde{A}_{k}\otimes \Big (F\otimes \psi ^*(E)^{\otimes
p^{k-1}}\Big ),
\end{eqnarray*}
where $\widetilde{A}_{k}=:\psi ^*(A)^{\otimes (p^{k-1}-k_0)}\otimes
\Big ( \psi ^*(A)^{k_0}(-F)\Big )$. By construction,
$\widetilde{A}_{k}$ is ample.

\begin{lem}\label{lem:new-estimate-power-of-p}
Fix $\varepsilon >0$. If $k\gg 0$, we have
$$\mathrm{vol}_\mu\big (\widetilde{A}_{k}\big )\ge
(p^{k})^{\mathrm{n}-1}\,\Big ( \mathrm{vol}_\mu(L)-\varepsilon \Big
).$$
\end{lem}

\textit{Proof.} By Lemma \ref{lem:exists-a-section}, for every $\mu
\in F$ we have
\begin{eqnarray*}
\mathrm{vol}_\mu\big (\widetilde{A}_{k}\big )&\ge
&\mathrm{vol}_\mu\big (A^{\otimes p^{k-1}-k_0}\big )\\
&=&\big (p^{k-1}-k_0\big ) ^{\mathrm{n}-1}\mathrm{vol}_\mu\big
(A_p\big )\\
&\ge&\big (p^{k-1}-k_0\big ) ^{\mathrm{n}-1}\,p^{\mathrm{n}-1}\,\Big
(
\mathrm{vol}_\mu(L)-\epsilon\Big )\\
&=&(p^{k})^{\mathrm{n}-1}\,\left (1-\frac{k_0}{p^{k-1}}\right
)^{\mathrm{n}-1}\,
\Big ( \mathrm{vol}_\mu(L)-\epsilon\Big )\\
&\ge&(p^{k})^{\mathrm{n}-1}\,\Big ( \mathrm{vol}_\mu(L)-\epsilon '
\Big ),
\end{eqnarray*}
where $\epsilon'$ can be made arbitrarily small by taking
$0<\epsilon \ll 1$ and $k\gg 0$.

\hfill Q.E.D.

This completes the proof of Proposition \ref{prop:fujita-step-2}.

\bigskip

From now on, we shall simplify our notation and set
$M'=\widehat{M}_\infty$, $\beta '=\widetilde{\beta}:M'\rightarrow
M$, $A'=\widetilde{A}_k$, and simply replace $p^k$ with $p$. Thus,
we have $(\beta')^*(L)^{\otimes p}=A'\otimes R$, for some
$G$-invariantly effective $R$, and $\mathrm{vol}_\mu(A')\ge
p^{\mathrm{n}-1}\,\big (\mathrm{vol}_\mu(L)-\epsilon\big )$ for
every $\mu \in F$.

\begin{prop}\label{prop:ample-stable-non-empty-on-resol}
The exists $B'\in \mathrm{Pic}^{S^1}(M')$ ample and with non-empty
stable locus: $(M')^\mathrm{s}\big (B'\big )\neq
\emptyset$.\end{prop}

\textit{Proof of Proposition
\ref{prop:ample-stable-non-empty-on-resol}.} Let $H\in
\mathrm{Pic}^{S^1}(M')$ be any ample Hermitian $S^1$-linearized line
bundle, with associated moment map $\Phi _H:M'\rightarrow
\mathbb{R}$  (by (\ref{eqn:moment-map-general}), the moment maps are
uniquely associated to the invariant Hermitian structure and the
linearization). Let $\Phi _B:M\rightarrow \mathbb{R}$ be a moment
map for $B$, with respect to an appropriate invariant Hermitian
structure. By the regularity of $B$, $\Phi _B^{-1}(0)\neq
\emptyset$, and $0$ is a regular value of $\Phi _B$. Let us define
$B_k=:\big (\beta '\big )^*(B)\otimes H\in \mathrm{Pic}^{S^1}(M')$.
Then $B_k$ has a product invariant Hermitian structure, with
associated moment map
\begin{equation}\label{eqn:moment-map-of-phi-k}
\Phi _{B_k}=k \cdot \big (\Phi _B\circ \beta'\big )+ \Phi _H= k\left
[\big (\Phi _B\circ \beta'\big )+ \frac 1k\Phi _H\right
].\end{equation}

Since $\beta'$ is a proper birational morphism, the image of its
exceptional locus is an invariant projective subvariety $Z\subseteq
M$. Thus $\beta'$ is an isomorpism over a dense open subset of $\Phi
^{-1}_B(0)$, so that we can find $m\in \Phi _B^{-1}(0)$ and an
analytic neighbourhood $M\supseteq U\ni m$ whose closure
$\overline{U}$ has positive distance from $Z$. Let us identify
$\overline{U}$ with its inverse image in $M'$. In view of
(\ref{eqn:moment-map-of-phi-k}), for $k$ sufficiently large
$\Phi_{B_k}$ is submersive at some point $m'\in
\Phi_{B_k}^{-1}(0)\cap \overline{U}\neq \emptyset$ near $m$. The
proof of Proposition \ref{prop:ample-stable-non-empty-on-resol} is
then completed by the following:

\begin{lem}\label{lem:regular-moment-map-implies-stable-not-empty}
Let $N$ be a projective $G$-manifold, and suppose $E\in
\mathrm{Pic}^G(N)$ is ample. Let $h=h_E$ be a $G$-invariant
Hermitian metric on $E$, and let $\Phi _E=\Phi_{A,h}$ be the induced
moment map (in the sense of (\ref{eqn:moment-map-general})). Suppose
that:
\begin{itemize}
    \item the normalized curvature form $\Omega$ of the associated compatible connection is
K\"{a}hler;
    \item $\Phi _E^{-1}(0)\neq \emptyset$, and
    $\Phi _E$ is submersive at some $m\in \Phi _E^{-1}(0)$.
\end{itemize}
Then $N^{\mathrm{s}}(E)\neq \emptyset$.
\end{lem}

\textit{Proof of Lemma
\ref{lem:regular-moment-map-implies-stable-not-empty}.} $h$ and
$\Omega$ induce unitary structures on $H^0(N,E^{\otimes k})$, $k\ge
0$. Let $\Omega _k$ be the associated Fubini-Study form on
$\mathbb{P}H^0(N,E^{\otimes k})^*$, and let $\varphi _k:M\rightarrow
\mathbb{P}H^0(N,E^{\otimes k})^*$ be the $G$-equivariant projective
embedding induced by the linear series $\big |E^{\otimes k}\big |$,
$k\gg 0$. Let $\Phi _k:\mathbb{P}H^0(N,E^{\otimes k})^*\rightarrow
\mathfrak{g}^*$ be the moment map for the $G$-action.

By a well-known theorem of Tian \cite{tian} and Zelditch \cite{z},
$\varphi _k$ is asymptotically isometric, in the sense that $\Omega
-\frac 1k\varphi _k^*\big (\Omega _k\big )=O\left (\frac{1}{k}\right
)$ in $\mathcal{C}^l$-norm, for every $l$. Therefore, if $k\gg 0$ we
also have $\Phi -\frac{1}{k}\Phi _k\circ \varphi _k=O\left
(\frac{1}{k}\right )$ in $\mathcal{C}^l$-norm, for every $l$.

It then follows that, for $k\gg 0$, $\big (\Phi _k\circ \varphi
_k\big ) ^{-1}(0)\neq \emptyset$, and that $\Phi _k\circ \varphi _k$
is submersive at some $m'\in \big (\Phi _k\circ \varphi _k\big )
^{-1}(0)$ close to $m$. Thus it suffices to prove statement is true
when $\Phi$ is a projectively induced moment map. This in turn
reduces the problem to proving it for a linear action on $
\mathbb{P}^N$, where the statement is a direct consequence of the
numerical criterion.

\hfill Q.E.D.

This completes the proof of Proposition
\ref{prop:ample-stable-non-empty-on-resol}.

\bigskip

To complete the proof of Theorem \ref{thm:equiv-fujita-approx}, we
now need only compose $\beta'$ with a Kirwan approximation, so as to
assume that there exist an ample regular line bundle on $M'$. That
the factorization can be lifted with the stated properties follows
by the same arguments as in the proof of Proposition
\ref{prop:fujita-step-2}.

\begin{cor} \label{cor:equiv-fuj-zero-vol}
Let $M$ be a projective $S^1$-manifold, and suppose that there
exists a regular and ample $B\in \mathrm{Pic}^{S^1}(M)$ such that
$\mathrm{codim}\big (M^\mathrm{u}(B)\big )\ge 2$. Let $\xi \in
\mathrm{NS}^{S^1}(M)$ be an $S^1$-big class and fix $\epsilon >0$.
Then there exists an equivariant modification $\beta :\tilde
M\rightarrow M$ (depending on $\xi$, $F$ and $\epsilon$) and classes
$a\in \overline{C^G(\widetilde{M})}$, $e\in
\mathrm{Eff}^{S^1}(\widetilde{M})$ such that $\beta^*(\xi)=a+e$ and
$ \mathrm{vol}_0(a)\ge \mathrm{vol}_0\big (\mu ^*(\xi)\big
)-\epsilon$; $a$ may be assumed to represent an ample line
bundle.\end{cor}

\end{document}